\documentclass[11pt]{article}

\usepackage{amsmath,amssymb}
\usepackage{hyperref}
\usepackage{graphicx}
\usepackage{float}
\usepackage{endnotes}
\usepackage{url}
\begin{document}

\title{ Mohsen Hachtroudi and the Geometry of Differential Equations}
\author{Masoud Khalkhali\footnote{\noindent Department of Mathematics, University of Western Ontario,  London Ontario, Canada. \\ masoud@uwo.ca}}
\date{}

\maketitle

\begin{abstract} 
Mohsen Hachtroudi (1908--1976) was one of the founders of modern
mathematics in Iran and  the first Iranian mathematician to gain
international recognition for original research in modern times.  A student of
\'Elie Cartan, he introduced the canonical projective connection now
known as the \emph{Hachtroudi connection}, whose significance has
become increasingly apparent through later developments in Cartan
geometry, CR geometry, and parabolic geometry.

 This article combines a historical and mathematical study of\\
Hachtroudi's life and work. It examines his education in Iran and
France, his role in establishing modern mathematical research and
higher education in Iran, and his influence as a teacher, essayist,
and public intellectual. It also presents a modern account of his
mathematical contributions and their continuing impact on the geometry
of differential equations and the theory of geometric structures.

 \end{abstract}

\tableofcontents

\section{Introduction}

Mohsen Hachtroudi (also transliterated as Hashtroudi or Hashtroodi) was
one of the founders of modern mathematics in Iran and  the first
Iranian mathematician to gain international recognition for original
research in modern times. His importance rests not only on his pioneering role in the
development of mathematics in Iran, but also on the lasting influence
of his mathematical contributions, which continue to resonate in
differential geometry, and the theory of differential equations.  Today the term \emph{Hachtroudi connection} appears naturally in the literature on Cartan geometry, the geometry of differential equations, projective differential geometry, CR geometry, and, more recently, parabolic geometries. Ironically, many geometers recognize the term  `Hachtroudi connection' without knowing much about the mathematician whose name it bears.

The aim of this article is twofold. On the one hand, it seeks to present an account of Hachtroudi's life and his decisive role in the development of modern mathematics in Iran. On the other hand, and perhaps more importantly, it aims to explain the mathematical significance of his work from a contemporary geometric viewpoint. Hachtroudi's doctoral thesis, written under the supervision of \'Elie Cartan in Paris in the 1930s, contains ideas that anticipated several themes which later became central in differential geometry. His construction of a canonical projective connection associated with completely integrable systems of second-order partial differential equations has become known as the \emph{Hachtroudi connection}.

Beginning with Shiing-Shen Chern's reinterpretation of Hachtroudi's
projective connection \cite{Chern75}, subsequent work by J\"urgen
Moser, Noboru Tanaka, David Grossman, and others revealed increasingly
deep connections between Hachtroudi's construction, Cartan's method of
equivalence, the geometry of differential equations, and CR and
parabolic geometries \cite{CM74,Tanaka76,Grossman00}. More recently,
Jo\"el Merker, partly in collaboration with Masoud Sabzevari, has
developed these connections further, both through the Cartan
equivalence problem for Levi-nondegenerate real hypersurfaces and
through the relation between CR geometry and systems of partial
differential equations \cite{MerkerLie,MerkerSabzevari12,
MerkerSabzevari14}.  In particular, Merker's work on the vanishing of Hachtroudi curvature
and local equivalence to the Heisenberg pseudosphere gives a direct
modern continuation of the geometric ideas originating in
Hachtroudi's thesis \cite{MerkerHach}.

Rather than presenting Hachtroudi's mathematics solely in its historical form, we shall interpret it from the perspective of modern differential geometry. This allows one to appreciate not only the originality of his ideas but also their continuing influence on contemporary mathematics.

This  paper is organized into two complementary parts. The first is
devoted to Hachtroudi's life and intellectual legacy: his education in
Iran and France, his return to Iran, his decisive role in the
development of modern mathematics and higher education in the country,
and his broader presence as a teacher, essayist, science promoter, and
public intellectual. The second part examines his mathematical work,
beginning with Cartan's method of equivalence and leading to
Hachtroudi's construction of a canonical projective connection and its
curvature. We then trace the later reappearance of these ideas in CR
geometry through the work of Chern and Moser and their subsequent
development in modern Cartan and parabolic geometry. \\

\noindent{\it Acknowledgement.}
This article is a long overdue and substantially expanded version of a lecture that I delivered at Azarbaijan University of Tarbiat Moallem in Tabriz, Iran, on the occasion of the centenary of the birth of Mohsen Hachtroudi. The lecture was presented at  the 37th Annual Iranian Mathematics Conference,  September 2–5, 2006.  I am grateful to the organizers of both events for the invitation and the opportunity to speak on Hachtroudi's mathematics and life. 
I should also    thank Siavash Shahshahani for  some  early discussions about   Hachtroudi's mathematical work,   and  Megerdich Toomanian for  providing me with a copy of the published version of Hachtroudi's Paris doctoral thesis.  
Finally, I wish to express my  gratitude to Hadi Soudbakhsh for
his  encouragement and support at an early stage, and for  sharing his
personal recollections of  Hachtroudi, both from his years as
one of his students at the University of Tehran and from their
subsequent close personal association.
 Needless to say, any historical or mathematical  inaccuracies  are entirely my own.

\section{Early Life and Education}

\begin{figure}[H]
   \centering
   \includegraphics[width=0.7\textwidth]{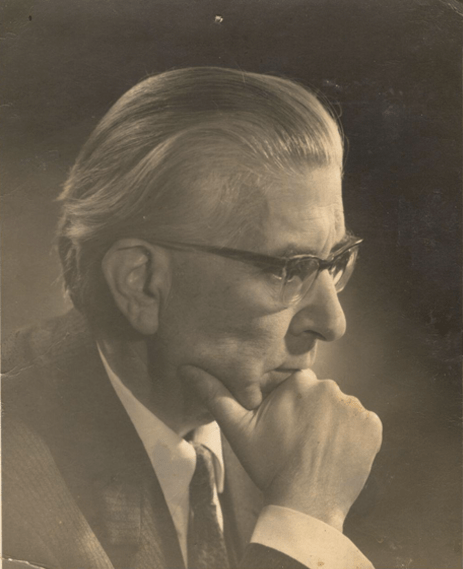} 
   \caption{Mohsen Hachtroudi (1908-1976)}
  \label{} \end{figure}

Mohsen Hachtroudi was born on  12 January  1908, in Tabriz, Iran,  and died on  4 September 1976, in Tehran.\footnote{Sources disagree on the exact conversion of Hachtroudi's birth date; the
Persian date 22 Dey 1286 corresponds to 12 January 1908, whereas several
secondary sources including the Wikipedia give 12 January 1907.}. He belonged to a distinguished family whose public and intellectual activities were closely connected with one of the most important periods of modern Iranian history. His father, Sheikh Esm\=a'il Mojtahed, was a prominent clergy and religious scholar and an adviser to Sheikh Mohammad Khi\=ab\=ani, one of the leading figures of the Constitutional Movement in Azerbaijan. The Constitutional Revolution (1905--1911) profoundly transformed the political and intellectual life of Iran, introducing parliamentary government and stimulating educational and cultural reforms. Hachtroudi thus grew up in an environment where scholarship, public service, and intellectual independence were deeply valued.

Equally influential was Hachtroudi's elder brother,
Mohammad Zia \\ Hachtroudi, whom Mohsen himself later described as the
person to whom he owed ``the foundation of my education''\cite{HachAuto}. Mohammad
Zia became a distinguished writer, editor, and literary critic, and
played an important role in the emergence of modern Persian
literature. He was among the earliest intellectuals to recognize the
importance of Nima Yushij and helped introduce his modern 
poetry to a wider audience through some of the earliest anthologies
of contemporary Persian literature.

It is a striking historical coincidence that the two Hachtroudi
brothers became pioneers of modernization in two different spheres of
Iranian intellectual life. While Mohammad Zia helped introduce modern
Persian literature to a new generation of readers, Mohsen became one
of the founders of modern mathematical research and education in
Iran. Both combined deep respect for the classical tradition with a
firm belief that intellectual renewal required engagement with the
most advanced ideas of their time.

Mohsen Hachtroudi's  early education took place in Tehran. After attending the Sirus and Aghdasieh primary schools, he entered the celebrated D\=ar al-Fon\=un, the first modern institution of higher learning in Iran. Founded in 1851 by Amir Kabir, D\=ar al-Fon\=un played a central role in introducing modern science and mathematics into Iran and educated many of the country's future scientists, engineers, physicians, and statesmen. Hachtroudi graduated from D\=ar al-Fon\=un in 1925, already displaying exceptional mathematical ability.

Hachtroudi's education abroad did not constitute a single uninterrupted period. After graduating from Dār al-Fonūn  he studied medicine for several years and then made his first trip to Europe, after which he returned to Iran. Upon his return, he enrolled in the Higher Teachers' College (\textit{Dār al-Mo`allemīn-e Markazī}), where he chose mathematics as his field of study and graduated as a member of its second graduating class. Only afterward did he make a second trip to France, where he studied at the Faculty of Sciences of the University of Paris. There he became a student of the great French geometer \'Elie Cartan (1869--1951), whose revolutionary work on Lie groups, differential geometry, moving frames, and the method of equivalence was transforming the subject.

\begin{figure}[H]
   \centering
   \includegraphics[width=0.6\textwidth]{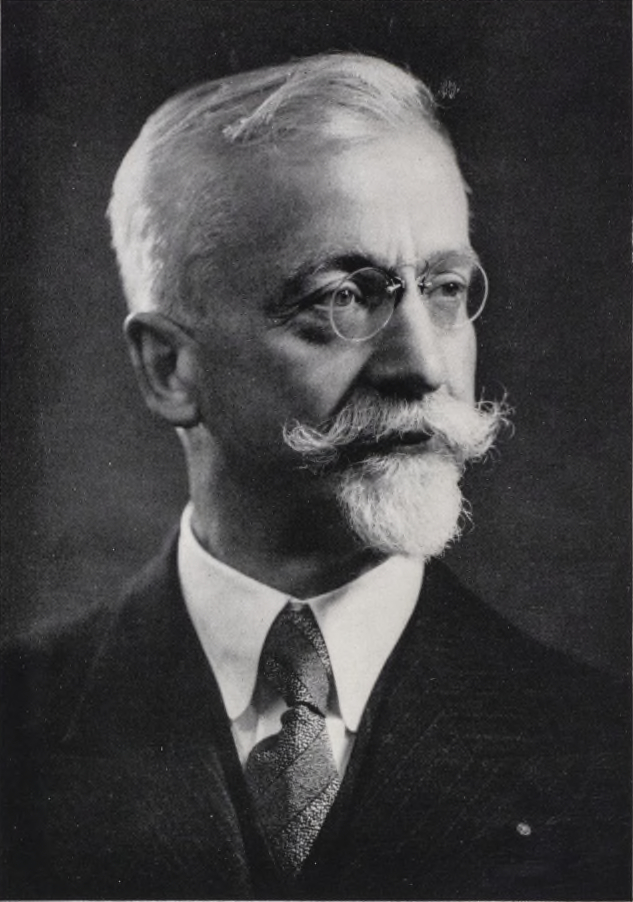} 
   \caption{ \'Elie Cartan }
  \label{} \end{figure}

In 1937 he was awarded the French \emph{Doctorat d'\'Etat},
then the highest academic doctorate in France, for his thesis
\begin{center}
\emph{Sur les espaces d'\'el\'ements \`a connexion projective normale}.
\end{center}
In this work he constructed the canonical projective connection that
would later become known as the \emph{Hachtroudi connection}.  His thesis  was published shortly afterwards in Paris in the same year\cite{Hach37}.

\begin{figure}[H]
   \centering
   \includegraphics[width=0.7\textwidth]{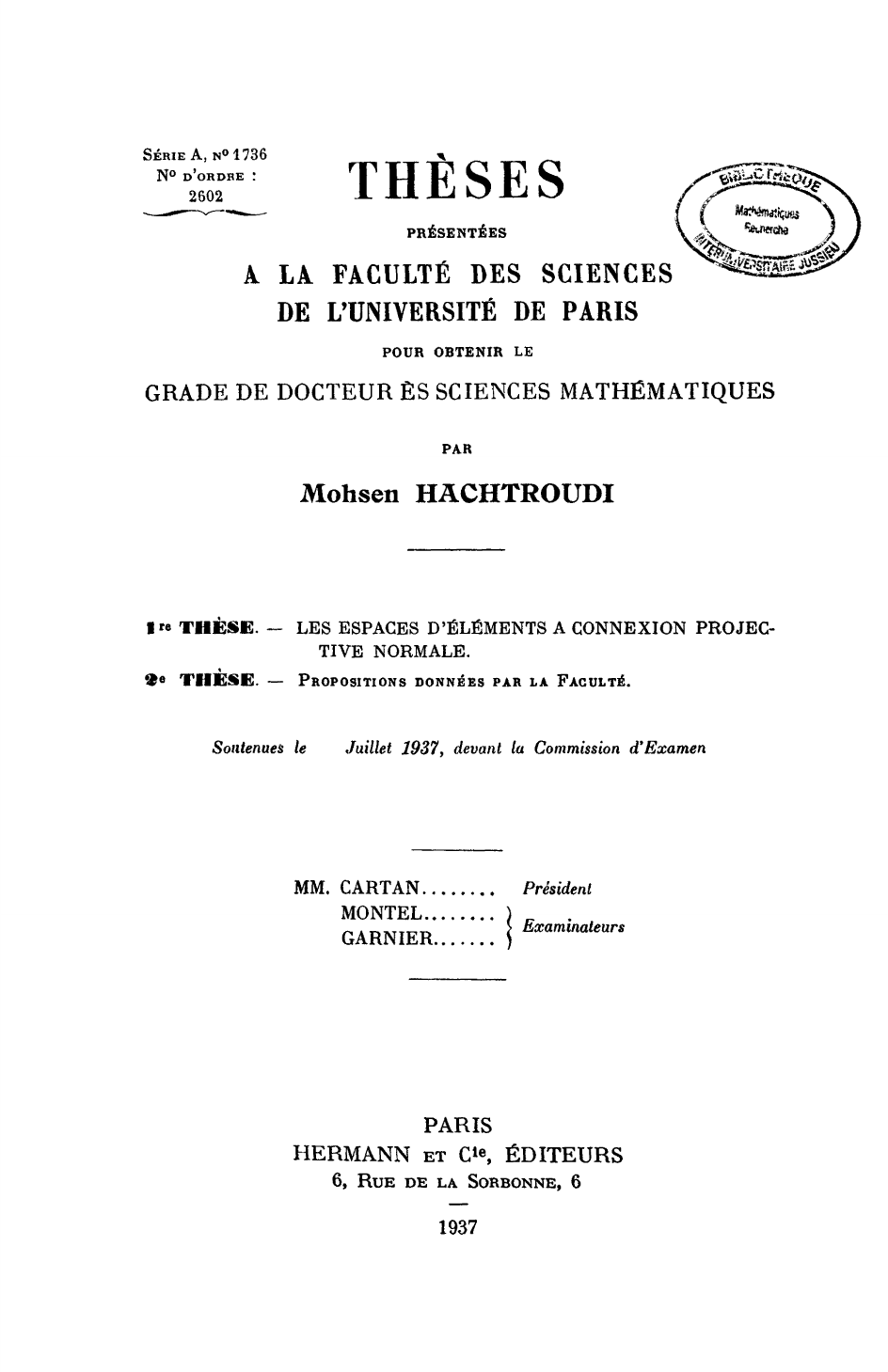} 
   \caption{The cover page of Hachtroudi's 1937 Paris thesis  written under \'Elie Cartan. The examiners were prominent French mathematicians of the time, Montel and Garnier. }
  \label{} \end{figure}

The thesis belongs to Cartan's program of attaching canonical geometric structures to differential equations and solving equivalence problems through curvature. Hachtroudi considered completely integrable systems of second-order partial differential equations and associated with them a canonical normal projective connection. The curvature of this connection provides the complete obstruction to transforming the system into its flat model. This construction, now commonly referred to as the \emph{Hachtroudi connection}, significantly generalized Cartan's earlier work\cite{Cartan24} and anticipated several ideas that became central decades later in the theories of CR geometry, Cartan geometries, and parabolic geometries.

\begin{figure}[H]
   \centering
   \includegraphics[angle=270, width=0.6\textwidth]{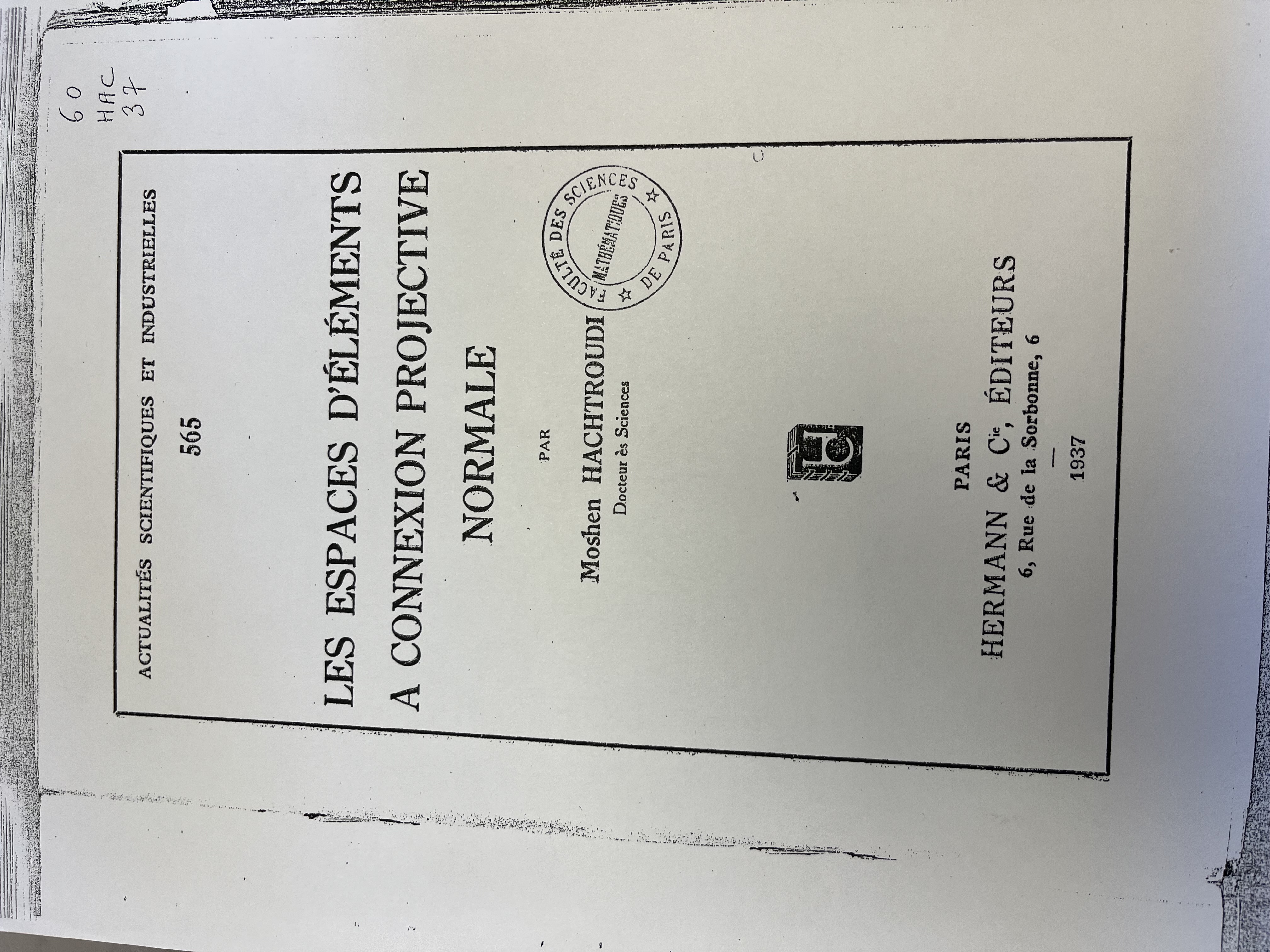} 
   \caption{The cover page of Hachtroudi's Paris thesis  published  by  Hermann in 1937.}
  \label{} \end{figure}

Although the immediate reception of Hachtroudi's thesis was limited, its influence steadily grew. In particular, Shiing-Shen Chern later recognized that the projective connection arising from the Segre family of a Levi nondegenerate CR hypersurface is, essentially, the same construction introduced by Hachtroudi. Through this observation, Hachtroudi's work became directly linked with the celebrated theory of Chern and Moser on the local equivalence of real hypersurfaces in complex manifolds.\\

\section{A Generation Shaped by Educational Reform}

Mohsen Hachtroudi belonged to a remarkable generation of Iranian
scholars whose education coincided with arguably  the most ambitious
programs of modernization in the country's history. During the
two decades between the rise of Reza Shah Pahlavi in 1925 and the end
of his reign in 1941, Iran underwent profound institutional, social,
and educational transformations. Pursuing ideals   of  the constitutional revolution of 1906  and its  founding fathers, the government of Reza Shah initiated  an extensive
program of state building that included the expansion of modern
administration, the construction of railways and roads, the
development of industry, and the establishment of new educational
institutions. Equally important was the emergence of a modern
professional middle class consisting of teachers, engineers,
physicians, lawyers, civil servants, and scientists, many of whom
would become the intellectual leaders of the country during the
following decades. 

Education occupied a central place in these reforms. Beginning in the
1920s, the government undertook a comprehensive reorganization of the
school system based largely on European, and especially French,
models. Modern elementary and secondary schools were established
throughout the country, teacher-training colleges were expanded, and
primary education was made compulsory. These reforms dramatically
increased school enrollment and created, for the first time, a
national system of secular public education reaching far beyond the
major cities. 

One of the most enduring achievements of this period in the area of education was the creation
of the University of Tehran. The law establishing the university was
approved by the Majles in March 1934, and the university was formally
inaugurated on 4 February 1935 with six faculties, including the
Faculty of Sciences, where Hachtroudi would later spend most of his
academic career. Conceived as the country's first modern university, the University of Tehran 
 brought together several existing professional schools into a
single institution devoted to both teaching and research. The
organization of the university followed largely the French model,
reflecting the educational ideals of many of its founders, including
Ali-Asghar Hekmat and Isa Sadiq, who had themselves studied in Europe.

A particularly important component of these educational reforms was
the establishment of a systematic program of state scholarships for
study abroad.  On 22 May 1928, the Iranian Parliament (Majlis)
passed the \emph{Law for Sending Students Abroad} (in Persian,
\textit{Q\=anun-e E'z\=am-e Mohassel be Kh\=arejeh}).
  The law required the government to
provide increasing annual appropriations for sending students from
both the capital and the provinces abroad for advanced study in
scientific and technical fields determined by the state.  Candidates
were to be selected through a competitive examination administered by
the Ministry of Education, with preference, under otherwise equal
conditions, given to holders of secondary, technical, or higher-school
diplomas.  The law also stipulated that at least thirty-five percent
of the students sent abroad each year were to study education and
teacher training.

Beginning with this program, talented young Iranians could compete for
state support to continue their education at leading European
universities, principally in France, Germany, Switzerland, and Great
Britain. Between 1928 and 1933
alone, approximately 640 students were sent abroad under this program.  The significance of the reform went beyond the number of
students who benefited from it.  By explicitly extending the program
to candidates from both Tehran and the provinces and making competitive
examination a central criterion of selection, it helped broaden access
to advanced European education beyond the relatively narrow privileged
circles for whom such opportunities had previously been available.
Hachtroudi belonged to the generation of young Iranian scholars whose
education and subsequent careers were profoundly shaped by this new
policy.

Mohsen Hachtroudi was one of the most distinguished members of this
generation. After completing his studies in Iran, he received a
government scholarship to continue his education in France, where he
entered the University of Paris and came under the supervision of
\'Elie Cartan. His subsequent career illustrates both the success of
this educational policy and the remarkable scientific achievements of
the first generation of Iranian mathematicians trained in the great
European mathematical schools.

The speed with which Hachtroudi mastered modern mathematics during
his years in Paris is all the more remarkable when viewed in the
context of mathematical education in Iran at the time. During the
1920s, higher education in mathematics was still in its infancy. The
curriculum at institutions such as D\=ar al-Fon\=un emphasized the
classical subjects required for engineering, military science, and
teacher training, while the revolutionary developments that had
transformed mathematics in Europe during the late nineteenth and early
twentieth centuries, concepts like Lie groups, differential geometry, topology,
functional analysis, and abstract algebra, were almost entirely
absent from university instruction. Research mathematics, in the
modern sense, had not yet taken root in Iran. 

Against this background, Hachtroudi's intellectual development was
nothing short of extraordinary. Within only a few years of study at
the University of Paris he not only acquired a thorough command of the
modern language of differential geometry but also entered one of the
most active research schools in the world under the direction of
\'Elie Cartan. His doctoral thesis demonstrates a complete mastery of
Cartan's method of equivalence, generalized spaces, projective
connections, and the geometry behind  differential equations. 

In fact the mathematical environment that Hachtroudi encountered in Paris was
exceptionally stimulating. During the 1930s the Faculty of Sciences of
the University of Paris was one of the world's foremost centers for
differential geometry. At its heart stood \'Elie Cartan, whose
revolutionary work on Lie groups, moving frames, generalized spaces,
and the method of equivalence had fundamentally reshaped modern
geometry. Around Cartan there had formed an outstanding school of
young geometers and analysts, including Charles Ehresmann, Andr\'e
Lichnerowicz, and later Shiing-Shen Chern, who would themselves become
leading figures in twentieth-century mathematics. Hachtroudi was thus
immersed in one of the most vibrant mathematical communities of the
period, where geometry, differential equations, topology, and
mathematical physics interacted in a remarkably fruitful way. 

Hachtroudi took full advantage of this unique intellectual
environment. In only a few years he progressed from a student educated
within the newly developing Iranian university system to a researcher
capable of contributing original ideas at the highest international
level. His doctoral dissertation exhibits not merely a mastery of
Cartan's methods but also an independence of thought that enabled him
to extend them in a substantial and original direction. It is perhaps
the clearest indication of his mathematical maturity that Cartan
himself later referred to Hachtroudi's thesis as ``une grande
généralisation'' of his own earlier work, while nearly forty years
later Shiing-Shen Chern identified Hachtroudi's projective connection
as the geometric structure underlying important aspects of the local
geometry of Levi-nondegenerate CR hypersurfaces. These two independent
assessments, coming from two of the twentieth century's greatest
geometers, leave little doubt that Hachtroudi's thesis occupies a
distinguished place in the development of Cartan's method of equivalence in modern differential
geometry. We shall return to this theme in more detail in the second part of this paper. 

\section{Return to Iran and Academic Career}

After receiving his \emph{Doctorat d'\'Etat} under \'Elie Cartan in
Paris in 1937, Hachtroudi returned to Iran at a decisive moment in the
country's educational history. The University of Tehran, founded only
a few years earlier in 1934 as the first modern university in Iran,
had become the centerpiece of an ambitious national program of higher
education and scientific modernization. Hachtroudi belonged to the
first generation of Iranian mathematicians trained at the highest
international level who returned home with the explicit goal of
building modern scientific research and university education in Iran.
He joined the Faculty of Sciences of the University of Tehran as an
Associate Professor and was promoted to full Professor in 1942.

Hachtroudi's academic career soon extended well beyond teaching and
research. In 1943 he was appointed Director of Education for Tehran,
reflecting the confidence placed in him by the country's educational
authorities. In 1952 he became President of the University of Tabriz,
where he played an important role in strengthening higher education in
northwestern Iran. Later, in 1958, he served as Dean of the Faculty of
Sciences of the University of Tehran. Throughout these years he
remained deeply committed to the development of mathematics, combining
research, teaching, and academic leadership with remarkable energy.

In 1945 Hachtroudi married Rob\=ab Modiri. They had three children:
two daughters, Far\=anak and Farib\=a, and one son, R\=amin.

\begin{figure}[H]
   \centering
   \includegraphics[ width=0.9\textwidth]{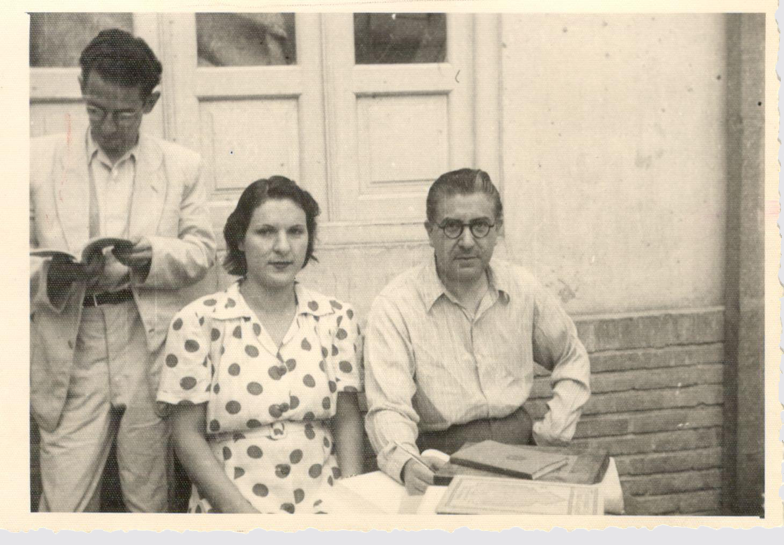} 
   \caption{At home with wife Rob\=ab Modiri and the Iranian writer Sadegh Hedayat. (Source: Association Mohsen Hachtroudi  \url{www.mo-ha.com})}
  \label{} \end{figure}

Scientifically, Hachtroudi continued to develop many of the ideas
initiated in his doctoral dissertation. His research ranged over
projective and affine differential geometry, Weyl geometry,
non-holonomic spaces, differential invariants of systems of
differential equations, analytical mechanics, Riccati equations, and
continued fractions. Much of this work remained firmly rooted in the
geometric philosophy of Cartan's school, seeking intrinsic geometric
structures underlying differential equations and transformation
groups. Several of his French monographs, published by the University
of Tehran, testify to the breadth and originality of these
investigations\cite{Hach45, Hach48, Hach56}.

Hachtroudi's concern with mathematical education was not confined to
the university.  He also wrote mathematics for younger readers and
high-school students, presenting elementary topics from a distinctly
modern and often surprisingly advanced point of view
\cite{Hach66,Hach70}.  These writings are characteristic of his broader
educational outlook: rather than treating school mathematics merely as
a collection of techniques, he sought to expose young readers to the
ideas, structures, and problems that give mathematics its intellectual
interest.

Although based primarily in Tehran, Hachtroudi maintained close
scientific contacts with mathematicians throughout Europe and North
America. His autobiography records periods of study abroad, including
a visit in  the autumn of 1951 as a member of the School of Mathematics
at the Institute for Advanced Study in Princeton\cite{IAS}. He also maintained scientific
correspondence and friendships with many leading mathematicians of his
time, among them Charles Ehresmann, Andr\'e Lichnerowicz, Jan
Schouten, Shiing-Shen Chern, Oscar Zariski, Dirk Struik, and several
others.\footnote{See Appendix A for an English translation of his short autobiographical essay.}  These international contacts enabled him to remain closely
connected with contemporary developments in geometry while helping to
introduce modern mathematical ideas into Iran. 

Perhaps Hachtroudi's greatest contribution was as a teacher and
builder of institutions. He regarded the university not merely as a
place for transmitting knowledge but as the foundation upon which a
modern scientific community could be built. Through his lectures,
supervision of students, and participation in academic life, he
introduced several generations of young Iranian mathematicians to
modern differential geometry and contemporary mathematical research.
Many of his students later became professors and researchers who
played leading roles in the development of mathematics at universities
throughout Iran. For this reason Hachtroudi is widely regarded as one
of the founders of modern mathematical research and education in
twentieth-century Iran.  

Beyond his scientific achievements, Hachtroudi was widely admired for
the breadth of his intellectual interests. He was deeply interested in
philosophy, literature, Persian poetry, the history of science, and
the cultural life of Iran. Through numerous public lectures, essays,
radio broadcasts, and articles in journals such as \emph{Yek\=an}, he
worked tirelessly to communicate modern scientific ideas to a broader
educated audience. This breadth of vision made him one of the leading
public intellectuals of his generation rather than simply a
distinguished mathematician. The Encyclopedia Iranica article  by 
Tahvildar-Zadeh and Majidi  gives  a concise  picture of Hachtroudi's social life,   and his activities and impact  beyond mathematics as a public intellectual and educator in Iran\cite{IranicaHach}. 

A particularly vivid portrait of Hachtroudi is given by Abbas Milani
in his biographical essay in \emph{Eminent Persians}
\cite{MilaniHashtrodi}. Milani portrays him not simply as an
outstanding mathematician, but as an unusually broad and charismatic
intellectual personality. He emphasizes Hachtroudi's extraordinary
mathematical gifts and memory, his independence of mind, his devotion
to his students, and the remarkable range of his interests in science,
philosophy, literature, and poetry. His daughter's memorable
observation that ``mathematics was his Mecca,'' quoted by Milani,
captures something essential about Hachtroudi's intellectual outlook:
mathematics was for him not merely a profession, but part of a larger
vision of nature and human existence. Milani's account thus provides
a valuable portrait of the man behind the mathematics and helps explain
the unusually deep impression Hachtroudi left on his students and on
the wider intellectual life of his generation.

Hachtroudi's  lasting  reception and influence in Iran  is reflected in the memorial volume edited by
his former student Hadi Soudbakhsh, which contains contributions by
many of Hachtroudi's students, colleagues, and friends, several of
whom later became leading mathematicians and professors at
universities throughout Iran \cite{Yad1}. Together with
four essays by Hachtroudi himself, the volume provides a vivid
portrait of his scientific achievements, his generosity as a teacher,
and his enduring influence on the intellectual and cultural life of
modern Iran.  This  influence  is particularly
evident in the special memorial issue of \emph{Yek\=an} published
shortly after his death, which brought together recollections and
assessments of his scientific, educational, and intellectual legacy
\cite{Yekan76}.

Hachtroudi's influence continued long after his death on 4 September
1976. In recognition of his scientific achievements and his decisive
role in the development of modern mathematics in Iran, the Iranian
Mathematical Society established the \emph{Professor Mohsen
Hachtroudi Prize} in his memory,  
awarded for distinguished work in geometry and topology.

\begin{figure}[H]
   \centering
   \includegraphics[angle=270,  width=0.7\textwidth]{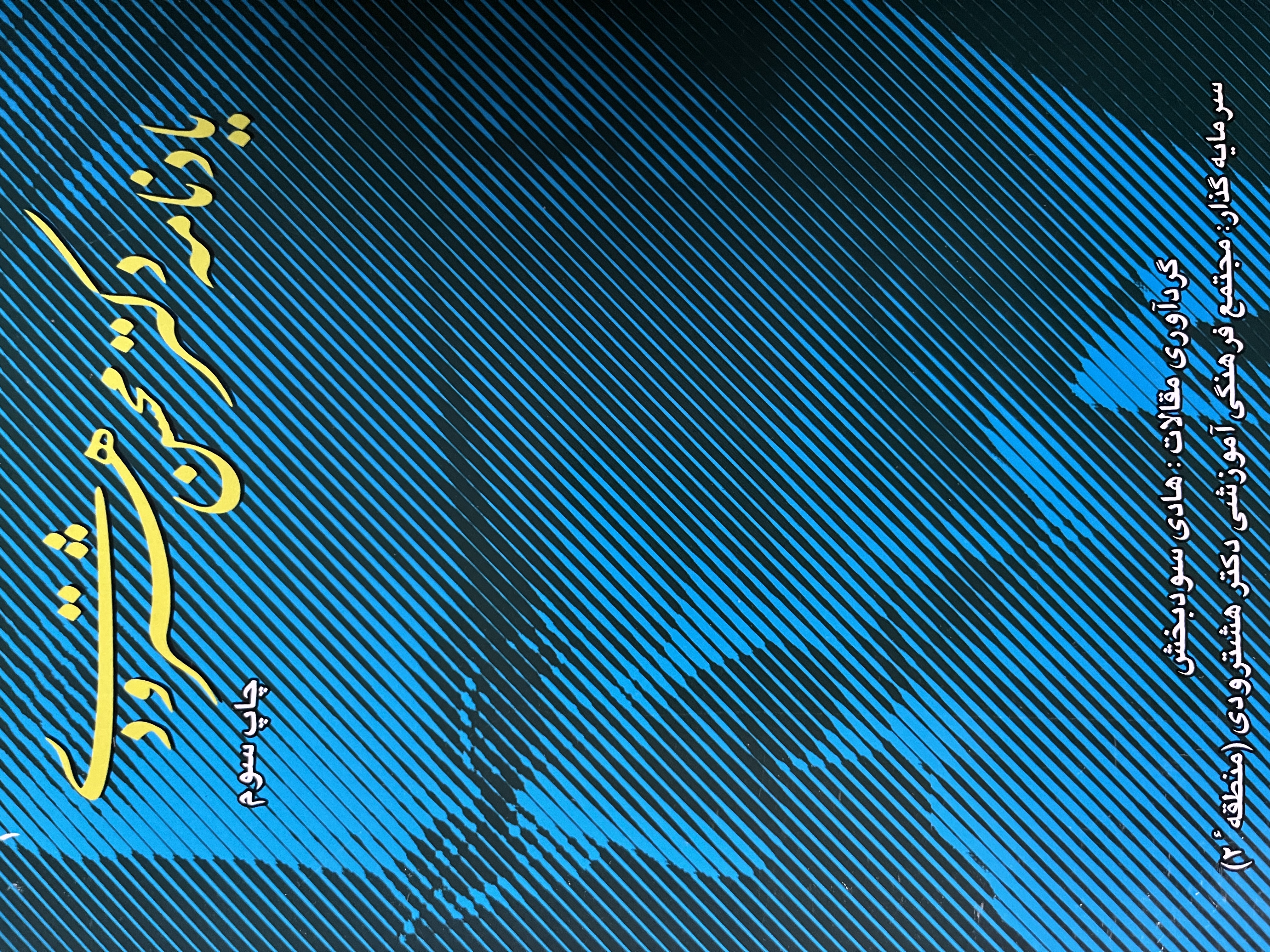} 
   \caption{The cover of the third edition (2007) of the memorial volume
\emph{Y\=adn\=ameh-ye Mohsen Hachtroudi}, edited by Hadi Soudbakhsh and
published in Tehran. The volume contains contributions by many of
Hachtroudi's former students, colleagues, and friends, several of whom
later became leading mathematicians and professors at Iranian
universities. It offers a unique and intimate portrait of Hachtroudi's
scientific, educational, and cultural legacy.}

\end{figure}

\section{Hachtroudi the Essayist}

To the general public in Iran, Mohsen Hachtroudi was known not only as
a distinguished mathematician but also as one of the country's leading
scientific essayists and public intellectuals. Throughout his career
he devoted considerable effort to communicating scientific ideas to a
broad audience. His lectures, essays, newspaper articles, and radio
broadcasts introduced generations of Iranian readers and listeners to
modern mathematics, physics, cosmology, philosophy of science, and the
relationship between science and culture.

Unlike many popularizers of science, Hachtroudi never regarded science
as an isolated intellectual activity. For him, scientific thought,
philosophy, literature, music, and the visual arts were different
manifestations of the same creative human spirit. This broad humanistic
vision characterizes much of his non-technical writing and explains why
his essays continue to be read long after the scientific questions that
motivated some of them have evolved.

The breadth of Hachtroudi's interests was not merely that of a
mathematician with literary tastes; mathematics, science, philosophy,
and poetry seem to have formed for him parts of a single intellectual
vision.  Milani aptly describes him as ``more in the spirit of a
Renaissance man than a monomaniacal mathematician'' \cite{MilaniHashtrodi}.

Perhaps the best illustration of this aspect of Hachtroudi's work is
the volume \emph{D\=ane\v{s} o Honar} (``Science and Art'') in Persian, first
published in 1961 and reprinted several times thereafter\cite{Hach61}. The book
collects a number of lectures and essays, some originally delivered
before public audiences and others previously published in newspapers
and magazines. Together they reveal the remarkable breadth of his
intellectual interests. Among the subjects discussed are the
relationship between science, literature, and art; the nature of
scientific and artistic creativity; the geometry and structure of the
universe; Einstein's theory of relativity and the nature of physical
reality; mathematical logic and the foundations of mathematics; the
principle of causality; the place of humanity in the cosmos and the
exploration of space; and contemporary Persian literature and literary
criticism. Read today, these essays reveal not only Hachtroudi's
scientific culture but also his remarkable ability to communicate
modern scientific ideas to a broad educated audience.

A  striking feature of the collection is its remarkable
range. Alongside essays on mathematical logic, relativity, and the
geometry of the universe, one finds thoughtful discussions of modern
Persian poetry, literary criticism, the search for extraterrestrial
life, the exploration of space, and even contemporary speculation
about flying saucers. This unusual combination reflects Hachtroudi's
conviction that science, philosophy, and the humanities form parts of
a single intellectual enterprise.

\begin{figure}[H]
   \centering
   \includegraphics[angle=270, width=0.8\textwidth]{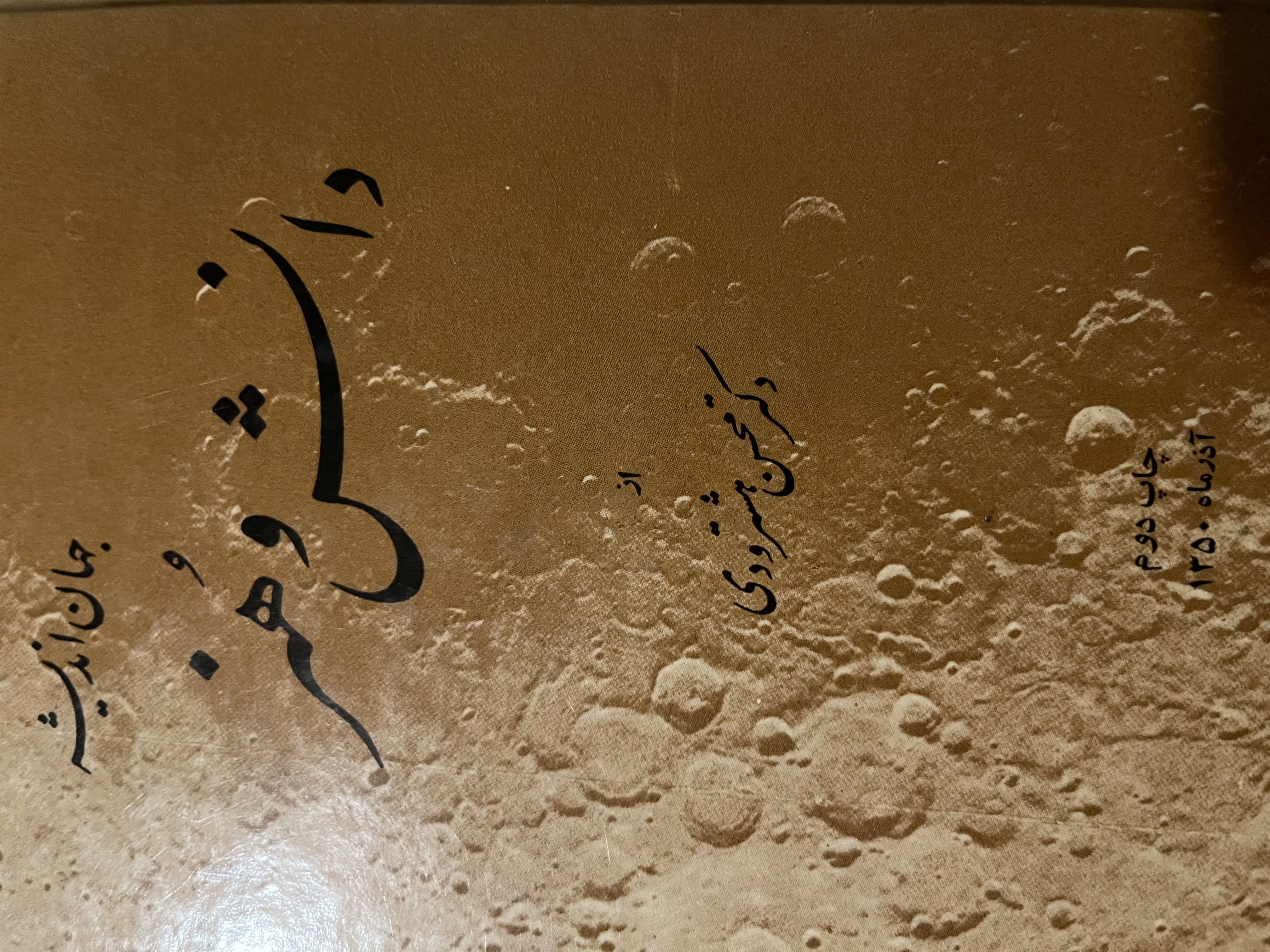} 
   \caption{The cover  of the second edition of Hachtroudi's book of essays  \emph{D\=ane\v{s} o Honar} (Science and Art) published in Iran first in 1961.}
  \label{} \end{figure}

Several themes recur throughout these essays. One is Hachtroudi's
conviction that scientific thinking is fundamentally universal, whereas
art reflects the particular historical and cultural experience of a
people. Another is his insistence that genuine scientific progress
requires creativity no less than artistic creation. Throughout the
volume he argues that science and art should not be viewed as opposing
activities but as complementary expressions of the human mind.

Another revealing example of Hachtroudi's broader intellectual
writing is his short book \emph{Seyr-e Andishe-ye Bashar}
(\emph{The Course of Human Thought}) \cite{Hach83}.  Unlike the
collection of essays in \emph{Science and Art}, this is a more focused
reflection on the historical development of scientific thought and on
the nature and teaching of science.  It illustrates particularly well
Hachtroudi's conviction that science should be understood not merely
as a body of established results, but as a continuing development of
human thought.

These ideas were not confined to his published essays. Hachtroudi was
also a gifted public speaker whose radio lectures introduced a wide
audience to contemporary developments in mathematics, physics, and
astronomy. Although only a small fraction of these broadcasts has
survived, they testify to his extraordinary ability to explain
difficult scientific ideas in clear, elegant Persian while preserving
their intellectual depth. His radio series \emph{Marzh\=a-ye D\=ane\v{s}}
(``Frontiers of Knowledge'') became particularly well known and is now
recognized as an important part of Iran's scientific and cultural
heritage.

One particularly revealing example of Hachtroudi's educational vision
appeared in the very first issue of \emph{Yek\=an}. Rather than
presenting elementary material, he chose to introduce his young
readers to the language of modern mathematics in an article entitled
``Foundations of Modern Mathematics.'' Beginning with the basic notions
of set theory, he proceeded to discuss the fundamental algebraic
structures of groups, rings, modules, and fields before introducing
the ideas of non-Euclidean geometry. Such a choice of topics was
extraordinary for its time. In the mid-1960s, these subjects were
still unfamiliar to many undergraduate mathematics students in Iran,
let alone high-school students. Hachtroudi firmly believed that
mathematically gifted students should be exposed, at an early stage,
to the conceptual foundations of modern mathematics rather than merely
to its computational techniques. In retrospect, this article was
remarkably ahead of its time and illustrates both his deep
understanding of contemporary mathematics and his visionary approach
to mathematical education.

\begin{figure}[H]
   \centering
   \includegraphics[ width=0.7\textwidth]{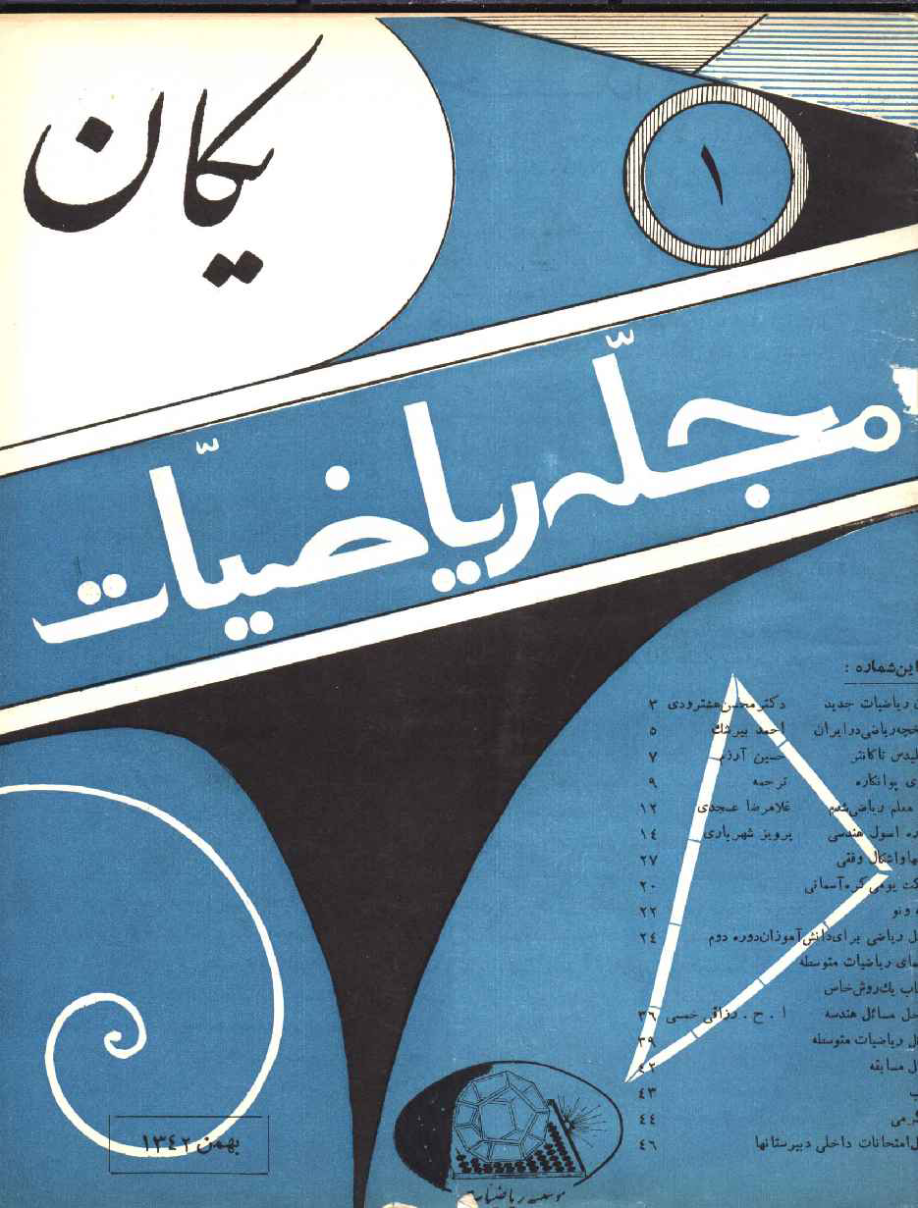} 
   \caption{\emph{Yek\=an}, one of Iran's most influential high-school mathematics magazines, published continuously from 1964 to 1977. Mohsen Hachtroudi was one of its principal contributors, publishing numerous expository articles, essays, and mathematical problems. Through \emph{Yek\=an}, he inspired several generations of young Iranian students and helped cultivate a lasting interest in mathematics and science. His very first article, appearing in the inaugural issue, was entitled ``Foundations of Modern Mathematics'' and introduced high-school students to such topics as set theory, groups, rings, modules, fields, and non-Euclidean geometry. At a time when many of these subjects were scarcely taught even in university mathematics programs in Iran, this was a remarkably forward-looking vision of mathematical education. (Source: Wikipedia)}
  \label{Yekan} \end{figure}

Hachtroudi's decision to begin with sets and algebraic structures also
reflects the international movement toward the modernization of
mathematics education that emerged during the 1950s and 1960s. His
article demonstrates that he was fully aware of these developments and
recognized the importance of introducing their central ideas to
Iranian students at an early stage.

Hachtroudi's influence therefore extended far beyond the university
classroom. For several decades he served as one of the principal
interpreters of modern science for the educated public in Iran. In this
respect, he belonged to a broader twentieth-century tradition
represented by figures such as Hermann Weyl, George Gamow, and Jacob
Bronowski: scientists and intellectuals who regarded reflection on the
meaning of science, and the communication of scientific ideas to a
wider public, as an essential part of their intellectual mission.

\section{A Bridge Between Two Mathematical Traditions}

The mathematical achievements of Mohsen Hachtroudi cannot be
understood apart from the remarkable historical moment in which he
lived. He belonged to the first generation of Iranian scholars whose
education coincided with the birth of modern higher education in Iran.
After receiving his scientific training under \'Elie Cartan in Paris,
he returned to devote his career to the development of mathematics in
his homeland. In doing so, he became a bridge between two great
mathematical traditions: the French school of modern differential
geometry and the rich but long-interrupted scientific tradition of
Iran, whose medieval mathematicians had once stood among the foremost
contributors to world mathematics. Hachtroudi did not merely import French differential geometry to Iran; he helped reconnect a great classical mathematical civilization with the modern international mathematical community

His scientific work, centered on the geometry of differential equations and the method of equivalence, continues to influence contemporary differential geometry nearly a century after its creation. What initially appeared as a highly specialized study of systems of second-order partial differential equations has become an integral part of several branches of modern geometry, including Cartan geometries, projective differential geometry, CR geometry, and parabolic geometries. Today the Hachtroudi connection is recognized as one of the earliest examples of a canonical Cartan geometry naturally associated with a nonlinear system of differential equations.

The subsequent sections of this article are devoted to explaining these mathematical developments in detail. Beginning with Cartan's method of equivalence, we shall describe Hachtroudi's construction of the canonical projective connection, its curvature and flatness criterion, its reinterpretation in terms of contact and incidence geometry, and finally its remarkable reappearance in CR geometry through the work of Chern and Moser. In this way we hope to demonstrate that Hachtroudi's work occupies a natural and enduring place in the history of modern differential geometry.

\section{From Cartan's Equivalence Problem to the \\Hachtroudi Connection}

Mohsen Hachtroudi's doctoral thesis represents one of the early
and remarkable applications of Cartan's method of equivalence to the
geometry of differential equations. Its principal construction, now
usually referred to as the \emph{Hachtroudi connection}, has subsequently
reappeared in projective differential geometry, CR geometry, the geometry
of differential equations, and, in modern terminology, the theory of
Cartan and parabolic geometries.

The purpose of this section is to explain the mathematical ideas behind
Hachtroudi's construction from a modern point of view before turning to
its more explicit formulation.

\subsection{Cartan and Chern on Hachtroudi's Thesis}

One of the strongest contemporary assessments of Hachtroudi's work
comes from his doctoral advisor, \'Elie Cartan himself. In his
\emph{Notice sur les travaux scientifiques}\cite{CartanNotice}, while discussing his own
work on projective connections associated with second-order ordinary
differential equations, Cartan concludes by observing that Hachtroudi
had recently obtained a far-reaching generalization of these results.
We reproduce below the relevant passage in the original French,
followed by an English translation.

\bigskip

\noindent

\begin{quote}\small

{\it Dans le mémoire [70]\footnote{Cartan's internal reference [70] corresponds to his 1924 paper \cite{Cartan24}}, j'ai considéré une équation différentielle
absolument arbitraire du second ordre
\[
y'' = f(x,y,y'),
\]
et je me suis demandé s'il est possible de regarder ses courbes
intégrales comme les droites d'un espace à deux dimensions à
connexion projective. Cela est possible d'une infinité de manières,
mais à la condition, si la fonction \(f\) n'est pas un polynôme du
troisième degré en \(y'\), de prendre comme élément générateur du
plan projectif, non le point, mais l'élément linéaire.

Si l'on veut maintenant trouver une connexion projective
intrinsèquement liée à l'équation donnée, j'ai démontré qu'il en
existe une, que j'ai encore appelée normale, jouissant de propriétés
géométriques remarquables, traduisant sous une forme géométrique les
propriétés analytiques de l'équation différentielle qui ne dépendent
pas du choix des variables.

On arrive tout naturellement à cette connexion, de même que dans les
deux exemples précédents, en appliquant ma méthode générale
d'équivalence à la recherche des invariants différentiels de
l'équation donnée par rapport au groupe infini des transformations
ponctuelles.}

\medskip

{\it Les résultats précédents ont reçu une grande généralisation dans
la thèse récente de M.~Hachtroudi
(\emph{Actualités Scientifiques}, Hermann, Paris, 1937).}

\end{quote}

\bigskip

\noindent
{\bf English translation:}

\begin{quote}

{\it In memoir [70], I considered an entirely arbitrary second-order
differential equation
\[
y'' = f(x,y,y'),
\]
and asked whether its integral curves could be regarded as the
straight lines of a two-dimensional space endowed with a projective
connection.

This is possible in infinitely many ways, provided that, whenever the
function \(f\) is not a polynomial of degree three in \(y'\), one
takes as the fundamental element of the projective plane not a point
but rather a line element.

If one now seeks a projective connection intrinsically associated
with the given differential equation, I showed that such a connection
indeed exists. I called it a \emph{normal} projective connection. It
possesses remarkable geometric properties which express, in geometric
form, the analytic properties of the differential equation that are
independent of the choice of variables.

This connection arises naturally, just as in the two preceding
examples, by applying my general method of equivalence to determine
the differential invariants of the given equation with respect to the
infinite group of point transformations.}

\medskip

{\it The preceding results have received a far-reaching
generalization in the recent doctoral thesis of
M.~Hachtroudi
(\emph{Actualités Scientifiques}, Hermann, Paris, 1937).}

\end{quote}

\medskip

From a modern perspective, this extension replaces the geometry of
path spaces by the geometry of jet spaces and incidence manifolds,
leading naturally to what is now known as the Hachtroudi connection.

Cartan's assessment was later reinforced by Shiing-Shen Chern, who placed Hachtroudi's thesis within the broader problem of associating a unique normal projective connection with a completely integrable family of submanifolds. In the introduction to the first volume of his \textit{Selected Papers}\cite{ChernSelectedI}, Chern described his own formulation of this problem and identified the case (k=n-1) as the main conclusion of Hachtroudi's thesis:
\begin{quote}

{\it The fundamental theorem on projective connections is the theorem associating
a unique normal projective connection to a system of paths. I announced in
[23] that the same is true when there is in a space of dimension \(n\) a family
of \(k\)-dimensional submanifolds depending on \((k+1)(n-k)\) parameters and
satisfying a completely integrable system of differential equations. The case
\(k=1\) is classical and the case \(k=n-1\) was the main conclusion of
M.~Hachtroudi's Paris thesis. My derivation was long and was never published.
A geometrical treatment was later given by C.~T.~Yen
(\textit{Annali di Matematica} 1953).}

\end{quote}

As was observed by Chern in his 1975 paper
\cite{Chern75}, the projective connection associated with the Segre
family of a Levi-nondegenerate real hypersurface is, apart from
notation, essentially Hachtroudi's normal projective connection.
Moreover, Chern showed that this projective connection underlies the
canonical connection introduced in the celebrated Chern--Moser theory
of CR hypersurfaces \cite{CM74}. This observation provides a direct
bridge between Hachtroudi's work on the geometry of systems of
differential equations and the modern differential geometry of CR
manifolds.  The equivalence problem studied by Hachtroudi belongs to a line of
investigation going back to the work of Lie and Tresse on differential
invariants of transformation groups \cite{LieEngel,Tresse}. Indeed,
Chern, in recalling Hachtroudi's construction, explicitly described it
as a generalization of the work of Tresse.

\subsection{Cartan's philosophy}

\'Elie Cartan (1869--1951) ranks among the greatest geometers of the
twentieth century and is widely regarded as one of the founders of
modern differential geometry.  For an in depth analysis of \'Elie Cartan's contributions to mathematics, see 
\cite{AkivisRosenfeld,ChernChevalley}. His work transformed several areas of
mathematics, including the classification and representation theory of
Lie groups and Lie algebras, the theory of symmetric spaces, Riemannian
and projective geometry, exterior differential systems, and the method
of moving frames. Perhaps his most enduring contribution was the
development of the method of equivalence, through which geometric
structures are characterized by canonical connections and their
curvature. These ideas profoundly influenced subsequent developments in
differential geometry, topology, mathematical physics, and the theory
of geometric structures. It was within this remarkable mathematical
environment that Hachtroudi carried out his doctoral research,
extending Cartan's ideas to completely integrable systems of
second-order partial differential equations and constructing the
canonical projective connection that now bears his name.

Cartan's method of moving frames provided a systematic way of replacing
an equivalence problem by the construction of an adapted coframe and
its associated structure equations \cite{Cartan35}.  This viewpoint,
closely connected with his theory of continuous transformation groups,
formed an essential part of the geometric framework in which
Hachtroudi's  thesis was conceived.

The familiar Levi--Civita connection of Riemannian geometry is one
example. Cartan's work showed that analogous canonical connections can
also be attached to projective and conformal geometries. Later, similar
ideas became central in CR geometry and in the general theory of Cartan
geometries.

The guiding principle may be expressed informally as follows: a geometric
structure should, whenever possible, be encoded by a canonical connection,
and its local equivalence problem should be translated into a question
about the curvature of that connection.
Thus, instead of producing a large collection of unrelated differential
invariants, one packages all the relevant local information into a single
geometric object.
Hachtroudi's thesis belongs naturally to this program.

The equivalence problem for differential equations has its origins in
Lie's theory of transformation groups \cite{LieEngel}, and was
subsequently transformed by Cartan into a powerful geometric method
based on moving frames, differential forms, and curvature.
For a modern treatment of differential equations from the viewpoint
of exterior differential systems and Cartan's geometric methods, see
\cite{BCGGG, Olver}.

It should also be mentioned that the geometry of differential equations had already become an active
subject in the 1920s, notably through the work of Veblen and Thomas on
the geometry of paths \cite{VeblenThomas23} and Cartan's theory of
projective connections. Hachtroudi's work belongs naturally to this
broader development of projective differential geometry.

\subsection{The equivalence problem for differential equations}

Hachtroudi considered systems of second-order partial differential
equations of the form
\[
\frac{\partial^2 y}{\partial x^i\partial x^j}
=
F_{ij}
\left(
x,y,
\frac{\partial y}{\partial x}
\right),
\qquad
1\leq i,j\leq n,
\]
where
\[
F_{ij}=F_{ji}.
\]
Here
$
x=(x^1,\ldots,x^n)
$
denotes the collection of independent variables, while
\[
y=y(x^1,\ldots,x^n)
\]
is the dependent variable. We shall often write this system as 
\[y_{ij}=F_{ij}(x,y,y_x).
\]
The equations are assumed to satisfy the appropriate compatibility
conditions so that the system is completely integrable.

The main problem is not to solve these equations explicitly. Rather, one
asks when two such systems represent the same local geometry.
Suppose another system is given by
\[
\frac{\partial^2 \widetilde y}
{\partial \widetilde x^i\partial \widetilde x^j}
=
\widetilde F_{ij}
\left(
\widetilde x,\widetilde y,
\frac{\partial \widetilde y}{\partial \widetilde x}
\right).
\]
One would like to determine whether there exists a local change of
variables
\[
(x,y)\longmapsto
(\widetilde x,\widetilde y)
\]
which transforms the first system into the second.
This is a typical example of what Cartan called an
\emph{equivalence problem}. The system of differential equations itself
is regarded as a geometric structure, and the objective is to find
quantities which remain invariant under the allowed transformations. The modern formulation of Cartan geometry provides a natural framework
for these constructions: a geometric structure is modeled on a
homogeneous space \(G/H\), while its deviation from the homogeneous
model is measured by the curvature of an associated Cartan connection;
see \cite{Sharpe}.

\subsection{The flat model}

As in many geometric equivalence problems, there is a distinguished
flat model. In the present case it is the system
\[
\frac{\partial^2 y}{\partial x^i\partial x^j}=0.
\]
Its solutions are the affine functions
\[
y(x)
=
a+\sum_{i=1}^{n}b_i x^i.
\]
Thus the solution hypersurfaces of the flat equation are affine
hyperplanes.
The fundamental question is therefore the following:
Under what conditions can a completely integrable system
\[
y_{ij}=F_{ij}(x,y,y_x)
\]
be transformed locally into
\[
y_{ij}=0
\]
by a change of variables?

This question is directly analogous to several classical questions in
geometry. For example, one may ask when a Riemannian metric is locally
Euclidean, when a projective structure is locally projectively flat, or
when a CR structure is locally equivalent to its homogeneous model.

In each of these situations, local flatness is characterized by the
vanishing of an appropriate curvature.
Hachtroudi discovered the corresponding curvature for the geometry of
the above systems of partial differential equations.

\subsection{Passing to first derivatives}

An important step is to enlarge the original space of variables
\[
(x^1,\ldots,x^n,y)
\]
by adjoining the first derivatives
\[
p_i=\frac{\partial y}{\partial x^i}.
\]
Thus one works on the space with coordinates
\[
(x^1,\ldots,x^n,y,p_1,\ldots,p_n).
\]
In modern terminology this is the first jet space
\[
J^1(\mathbb{R}^n,\mathbb{R}).
\]
The differential equation then becomes a geometric structure on this
jet space.
Introduce the contact form
\[
\theta
=
dy-\sum_{i=1}^{n}p_i\,dx^i.
\]
Along the first jet of a genuine function \(y=y(x)\), one has
\[
\theta=0.
\]
The second-order system determines the differential forms
\[
\theta_i
=
dp_i-\sum_{j=1}^{n}F_{ij}\,dx^j.
\]
Along a solution of the differential equation,
\[
\theta=0,
\qquad
\theta_i=0.
\]
Equivalently, the system determines the vector fields
\[
D_i
=
\frac{\partial}{\partial x^i}
+
p_i\frac{\partial}{\partial y}
+
\sum_{j=1}^{n}
F_{ij}
\frac{\partial}{\partial p_j}.
\]
These are the total derivative operators associated with the system.

Now complete integrability may be expressed geometrically by requiring the
distribution generated by
\[
D_1,\ldots,D_n
\]
to be integrable. In the simplest formulation this amounts to
\[
[D_i,D_j]=0,
\]
or to the corresponding compatibility conditions on the functions
\(F_{ij}\).

Thus a second-order PDE system has been converted into a geometric
distribution on a first jet space.

\subsection{From differential equations to projective geometry}

The crucial insight in Hachtroudi's work is that this geometric structure
carries a natural projective connection. The adjective `projective' here refers to the homogeneous
point--hyperplane incidence geometry
\[
G/P_{1,n+1},
\qquad
G=\mathrm{PGL}(n+2),
\]
rather than to an ordinary projective class of affine connections on the jet
space. The two transverse Legendrian foliations deform the two natural
fibrations of the projective incidence variety, and Hachtroudi's connection is
the canonical Cartan connection associated with this curved incidence
geometry.

The integral manifolds associated with the differential equation may be
viewed, in an appropriate geometric formulation, as playing a role
analogous to distinguished geodesic submanifolds. This suggests that the
correct invariant geometry is not primarily affine but projective.

Hachtroudi showed that one can associate to the differential equation a
canonical normal projective connection. The coefficients of this
connection are constructed from the functions \(F_{ij}\) and their
derivatives.
His construction is canonical: once the differential equation is given,
there are no arbitrary choices remaining in the resulting normal
connection.
It is this canonical projective connection which is now called the
\emph{Hachtroudi connection}.

\subsection{Curvature and the obstruction to flatness}

The main advantage of introducing a canonical connection is that the
equivalence problem can then be formulated in terms of curvature.
Let us recall that for  a Cartan connection \(\omega\), its curvature is
\[
\Omega
=
d\omega+ \frac{1}{2}\, \omega\wedge\omega.
\]
The flat model corresponds to
\[
\Omega=0 \quad \quad \text{on an open neighbourhood.}
\]
  Consequently, the curvature of the Hachtroudi connection measures the
failure of the original differential equation to be locally equivalent
to the flat system
\[
y_{ij}=0.
\]
Schematically, one has
\[
\text{Hachtroudi curvature}=0 \quad \text{locally}
\quad\Longleftrightarrow\quad
\text{local equivalence to the flat model}.
\]
This is the central geometric content of Hachtroudi's construction.\footnote{It should be stressed that vanishing of the curvature form $\Omega$ at a single point does not imply local flatness at all. So in the above statement and in similar statements in this paper,  by vanishing of curvature we mean vanishing locally in an open neighbourhood.} 

As we said before, the analogy with Riemannian geometry is particularly instructive. There,
the Riemann curvature tensor measures the obstruction to finding local
coordinates in which the metric becomes Euclidean. Here, Hachtroudi's
curvature measures the obstruction to finding local coordinates in which
the differential equation becomes
$y_{ij}=0.$

The conceptual achievement of Hachtroudi's work is therefore not merely
the computation of a particular invariant. It is the recognition that an
equivalence problem for a nonlinear system of partial differential equations can be encoded
in the curvature of a canonical geometric connection.
This construction may be summarized as
\[
\text{Second-order PDE}
\longrightarrow
\text{Geometry on a jet space}
\longrightarrow\]
\[\text{Canonical projective connection}
\longrightarrow
\text{Curvature invariants}.
\]
This is precisely the kind of geometric passage that later became
standard in the modern theory of differential equations.
It also explains why Hachtroudi's work has resurfaced in several areas of
contemporary geometry. The same general philosophy appears in the study
of CR structures, path geometries, Cartan geometries, parabolic
geometries, and overdetermined systems of partial differential equations.

From a modern perspective, Hachtroudi's thesis may therefore be viewed as
an early and important example of the principle that differential
equations themselves carry intrinsic geometry.

\section{Integrability, the Two Foliations, and the Flat Model}

We now examine more closely the geometric structure carried by the
first jet space. This is the natural setting in which the Hachtroudi
connection begins to emerge.
Throughout this section we consider a completely determined system
of second-order equations
\[
y_{ij}
=
F_{ij}(x,y,p),
\qquad
p_i=\frac{\partial y}{\partial x_i},
\qquad
1\leq i,j\leq n,
\]
where
\[
F_{ij}=F_{ji}.
\]
The underlying first jet space has coordinates
\[
(x^1,\ldots,x^n,y,p_1,\ldots,p_n)
\]
and dimension
$2n+1.$

Associated with the system are the total derivative operators
\[
D_i
=
\frac{\partial}{\partial x^i}
+
p_i\frac{\partial}{\partial y}
+
\sum_{k=1}^{n}
F_{ik}\frac{\partial}{\partial p_k},
\qquad
1\leq i\leq n.
\]
These vector fields have a simple meaning. If
\[
y=y(x)
\]
is a solution, then its first jet
\[
j^1y(x)
=
\left(
x,y(x),
\frac{\partial y}{\partial x^1},
\ldots,
\frac{\partial y}{\partial x^n}
\right)
\]
is tangent to the distribution
\[
E
=
\operatorname{span}\{D_1,\ldots,D_n\}.
\]
Thus the original differential equation is encoded geometrically by
an \(n\)-dimensional distribution \(E\) on the \((2n+1)\)-dimensional
first jet space.

Now the first question is whether this distribution is integrable and we shall address this question next. 

\subsection{The compatibility equations}

A direct computation gives
\[
[D_i,D_j]
=
\sum_{k=1}^{n}
\left(
D_iF_{jk}-D_jF_{ik}
\right)
\frac{\partial}{\partial p_k}.
\]
Consequently,
\[
[D_i,D_j]=0
\]
if and only if
\[
D_iF_{jk}
=
D_jF_{ik}
\]
for every \(i,j,k\).
Equivalently, after relabeling the indices, the compatibility
conditions may be written
\[
D_kF_{ij}
=
D_jF_{ik}.
\]
These are precisely the conditions obtained by differentiating the
original system and requiring equality of third derivatives. Indeed,
from
\[
y_{ij}=F_{ij}(x,y,y_x)
\]
one obtains
\[
y_{ijk}=D_kF_{ij},
\]
and the equality
\[
y_{ijk}=y_{ikj}
\]
therefore gives
\[
D_kF_{ij}=D_jF_{ik}.
\]
Thus complete integrability of the differential equation is equivalent
to the Frobenius condition $
[D_i,D_j]=0
$  derived above. 

\subsection{The foliation by solutions}

By the Frobenius theorem, the distribution
\[
E=\operatorname{span}\{D_1,\ldots,D_n\}
\]
is locally tangent to a foliation by \(n\)-dimensional submanifolds.
Each such leaf is the first jet of a solution of the differential
equation.
Locally the general solution can therefore be written in the form
\[
y
=
Q(x^1,\ldots,x^n,a^1,\ldots,a^n,b),
\]
where
\[
(a^1,\ldots,a^n,b)
\]
are \(n+1\) integration constants.
On the corresponding leaf one has
\[
p_i
=
\frac{\partial Q}{\partial x^i},
\]
and hence
\[
\frac{\partial^2 Q}{\partial x^i\partial x^j}
=
F_{ij}
\left(
x,Q,Q_{x^1},\ldots,Q_{x^n}
\right).
\]
Thus the leaves of the Frobenius foliation may be represented as
\[
x
\longmapsto
\left(
x,
Q(x,a,b),
Q_{x^1}(x,a,b),
\ldots,
Q_{x^n}(x,a,b)
\right).
\]
The parameter space
\[
(a^1,\ldots,a^n,b)
\]
therefore has dimension \(n+1\) and may be regarded locally as the
space of solutions.
For the flat equation
\[
y_{ij}=0,
\]
the general solution is
\[
Q(x,a,b)
=
b+\sum_{i=1}^{n}a^i x^i.
\]
Thus the solution space of the flat equation is simply the space of
affine hyperplanes in an \((n+1)\)-dimensional affine space.
This elementary observation is the first indication that projective
geometry is hidden in the problem.

\subsection{The contact structure}

The first jet space carries a canonical contact form
\[
\theta
=
dy-\sum_{i=1}^{n}p_i\,dx^i.
\]
Its kernel
\[
H=\ker\theta
\]
is a \(2n\)-dimensional contact distribution.
Since
\[
\theta(D_i)=0,
\]
the distribution \(E\) determined by the differential equation lies
inside the contact distribution:
\[
E\subset H.
\]

There is also a second, completely canonical \(n\)-dimensional
distribution
\[
V
=
\operatorname{span}
\left\{
\frac{\partial}{\partial p_1},
\ldots,
\frac{\partial}{\partial p_n}
\right\}.
\]
This is the vertical distribution of the projection
\[
J^1
\longrightarrow
\mathbb{R}^{n+1},
\qquad
(x,y,p)\longmapsto(x,y).
\]
Clearly,
\[
V\subset H.
\]
Moreover,
\[
E\cap V=\{0\},
\]
and, since both spaces have dimension \(n\),
\[
H=E\oplus V.
\]
Thus the contact distribution possesses a distinguished splitting into
two transverse \(n\)-plane distributions.

\subsection{The two distributions are Legendrian}

The exterior derivative of the contact form is
\[
d\theta
=
-\sum_{i=1}^{n}dp_i\wedge dx^i.
\]
Restricted to the contact distribution \(H\), this defines a
nondegenerate symplectic form, up to multiplication by a nonzero
function.
The vertical distribution \(V\) is isotropic since
\[
d\theta
\left(
\frac{\partial}{\partial p_i},
\frac{\partial}{\partial p_j}
\right)
=
0.
\]
Since \(\dim V=n\), it is maximally isotropic and hence Legendrian.
The same is true for \(E\). Indeed,
\[
d\theta(D_i,D_j)
=
-F_{ji}+F_{ij}
=
0,\]
because $
F_{ij}=F_{ji}.
$
Since $
\dim E=n,
$
the distribution \(E\) is also Legendrian.
We therefore obtain a splitting
\[
H=E\oplus V
\]
of the contact distribution into two transverse Legendrian
distributions.
When the system is completely integrable, both distributions are
integrable. The vertical distribution \(V\) is automatically
integrable, while the integrability of \(E\) is exactly the
compatibility condition discussed above.

In modern terminology, this structure is often described as a
bi-Legendrian structure, or, in a closely related formulation, as a
{\it generalized path geometry}.

\subsection{The two foliations}

The geometry therefore contains two distinguished foliations.
The first is the foliation $
\mathcal{F}_{\mathrm{sol}}
$
tangent to \(E\). Its leaves are the lifted solutions of the
differential equation.
The second is the vertical foliation $
\mathcal{F}_{\mathrm{vert}}
$ tangent to \(V\). Its leaves are obtained by fixing $
(x,y) $
and varying $
(p_1,\ldots,p_n).
$

Thus through each point of the jet space pass two distinguished
\(n\)-dimensional submanifolds, one belonging to each foliation.
This double-foliation structure is one of the essential geometric
features of Hachtroudi's problem.
Schematically,
\[
\begin{array}{ccc}
& J^1 & \\
\swarrow & & \searrow \\[2mm]
\text{space of points }(x,y)
& &
\text{space of solutions }(a,b).
\end{array}
\]
The left projection forgets the first derivatives \(p_i\), while the
right projection sends a jet to the unique local solution leaf
containing it.
Thus the jet space plays the role of an incidence space between
points and solutions.

This picture becomes especially transparent for the flat system $
y_{ij}=0, $ whose 
solutions are
\[
y=b+\sum_{i=1}^{n}a^i x^i.
\]
A point $
(x,y) $
belongs to the solution determined by \((a,b)\) precisely when $
y=b+\sum_{i=1}^{n}a^i x^i.
$
Thus the flat geometry is the incidence geometry between points and
affine hyperplanes.

After projective compactification, a point of affine
\((n+1)\)-space becomes a point of projective space $
\mathbb{P}^{n+1}, $
while an affine hyperplane becomes a projective hyperplane.
Consequently, the flat model is naturally identified with the
incidence variety
\[
\mathcal{I}
=
\left\{
(\ell,H):
\ell\subset H
\right\},
\]
where \(\ell\) is a line in \(\mathbb{R}^{n+2}\) and \(H\) is a
hyperplane containing it.
Equivalently,
\[
\mathcal{I}
=
\left\{
([X],[\Xi])
\in
\mathbb{P}^{n+1}\times(\mathbb{P}^{n+1})^*
:
\Xi(X)=0
\right\}.
\]
Its dimension is $
2n+1, $
exactly the dimension of the first jet space.

The two natural projections
\[
\mathcal{I}\longrightarrow\mathbb{P}^{n+1}
\]
and
\[
\mathcal{I}\longrightarrow(\mathbb{P}^{n+1})^*
\]
give precisely the two foliations appearing above.

The first fixes a point and varies the hyperplanes passing through
it. The second fixes a hyperplane and varies the points lying on it.
Thus the two foliations of the jet space for  nonlinear deformations
of the elementary point-hyperplane incidence geometry of projective
space are dual to each other. 

\subsection{The projective symmetry group}

The incidence relation
\[
\Xi(X)=0
\]
is preserved by projective transformations.
The natural symmetry group of the flat model is therefore
\[
PGL(n+2,\mathbb{R}),
\]
or locally, after passing to the corresponding Lie algebra,
\[
\mathfrak{sl}(n+2,\mathbb{R}).
\]
In the complex analytic setting one replaces these by
\[
PGL(n+2,\mathbb{C})
\]
and
\[
\mathfrak{sl}(n+2,\mathbb{C}).
\]
The incidence variety is a homogeneous space for this group. In
modern notation it may be written as 
\[
\mathcal{I}
\simeq
G/P,
\]
where
\[
G=PGL(n+2)
\]
and \(P\) is the parabolic subgroup preserving a flag consisting of a
line contained in a hyperplane.
Equivalently, using \(SL(n+2)\), the relevant parabolic is often
denoted
\[
P_{1,n+1}.
\]
This observation explains the appearance of projective geometry in
Hachtroudi's construction.
The flat differential equation $
y_{ij}=0 $
is not merely a particularly simple PDE. It is the affine realization
of a homogeneous projective incidence geometry.

\subsection{From the flat model to the Hachtroudi connection}

We can now see the basic idea behind Hachtroudi's construction.
For a general completely integrable system
\[
y_{ij}=F_{ij}(x,y,p),
\]
the first jet space still carries

\begin{itemize}
\item a contact distribution \(H\),
\item a decomposition
\[
H=E\oplus V,
\]
\item an integrable Legendrian distribution \(E\) whose leaves are
solutions,
\item an integrable Legendrian distribution \(V\) whose leaves are
vertical fibers.
\end{itemize}
For the flat equation, precisely this structure is realized by the
homogeneous projective incidence space.

The natural question is therefore whether the curved structure
associated with an arbitrary system can be compared canonically with
the flat incidence geometry.
Cartan's general principle suggests the answer: construct a principal
bundle
\[
\mathcal{G}\longrightarrow J^1
\]
together with a Lie-algebra-valued one-form
\[
\omega
\in
\Omega^1
\left(
\mathcal{G},
\mathfrak{sl}(n+2)
\right)
\]
which reproduces, infinitesimally, the geometry of the homogeneous
model.
The form \(\omega\) is a Cartan connection of type
\[
\bigl(\mathrm{SL}(n+2,\mathbb{F}),\,P_{1,n+1}\bigr),
\qquad
\mathbb{F}=\mathbb{R}\ \text{or}\ \mathbb{C}.
\]
Hachtroudi's achievement was to show that the differential equation
determines such a connection canonically, after imposing the
appropriate normalization conditions.
This is the \emph{Hachtroudi connection}.
Its curvature
\[
\Omega
=
d\omega+\frac{1}{2}\, \omega\wedge\omega
\]
measures the failure of the curved double-foliation geometry to be
locally equivalent to the flat incidence geometry.
Thus, as mentioned before,  one obtains the fundamental implication
\[
\Omega=0 \quad \text{locally}
\quad\Longleftrightarrow\quad
y_{ij}=F_{ij}
\ \text{is locally equivalent to}\
y_{ij}=0.
\]
This gives a geometric solution to the original equivalence problem.

\subsection{A useful dictionary}

It is useful to summarize the correspondence between the differential
equation and the associated geometry:

\[
\begin{array}{r@{\qquad\Longleftrightarrow\qquad}l}

\multicolumn{1}{c}{\textbf{Differential Equations}}
&
\multicolumn{1}{c}{\textbf{Geometry}}
\\[1.5ex]

y_i=p_i
&
\text{first jet coordinates}
\\[1ex]

y_{ij}=F_{ij}
&
\text{distribution }E
\\[1ex]

F_{ij}=F_{ji}
&
E\text{ is a Legendrian distribution}
\\[1ex]

D_kF_{ij}=D_jF_{ik}
&
E\text{ is integrable}
\\[1ex]

\text{solutions }y=Q(x,a,b)
&
\text{leaves of }E
\\[1ex]

(x,y)\text{ fixed}
&
\text{leaves of }V
\\[1ex]

H=E\oplus V
&
\text{double Legendrian structure}
\\[1ex]

y_{ij}=0
&
\text{flat path geometry}
\\[1ex]

\text{general }F_{ij}
&
\text{curved path geometry}
\\[1ex]

\text{Hachtroudi connection}
&
\text{canonical Cartan connection}
\\[1ex]

\text{Hachtroudi curvature}
&
\text{obstruction to local flatness}

\end{array}
\]

This dictionary contains, in compressed form, the main conceptual
content of Hachtroudi's construction. The remaining task is to make
the last two entries explicit: one must construct the Cartan
connection, impose its normalization conditions, compute its
curvature, and identify the curvature components which obstruct
equivalence with the flat equation. This will be our next step. The local equivalence problem for  completely integrable systems
of second-order partial differential equations has subsequently been
revisited from a modern geometric viewpoint; see, for example,
\cite{Bieche}.  Modern formulations of the equivalence problem for systems of
second-order differential equations were further developed by Fels
\cite{Fels95}, while Grossman's work on torsion-free path geometries
and integrable second-order systems placed closely related structures
within a broader geometric framework \cite{Grossman00}.

\section{From Hachtroudi's Connection to CR Geometry}

One of the most striking later developments of Hachtroudi's work is
its appearance in CR geometry as observed first by Chern in 1975\cite{Chern75}. At first sight, the geometry of
completely integrable systems
\[
y_{ij}=F_{ij}(x,y,y_x)
\]
and the geometry of real hypersurfaces in complex space seem to be
rather different subjects. The connection between them is provided by
the family of Segre varieties associated with a real analytic
hypersurface.

This relationship makes it possible to pass naturally from a CR
hypersurface to a system of second-order differential equations of
exactly the type studied by Hachtroudi. The Hachtroudi curvature of
the resulting system then becomes, in another language, the curvature
measuring the failure of the CR hypersurface to be locally equivalent
to a sphere or, more generally, to a nondegenerate hyperquadric.

This point of view also clarifies the relation between Hachtroudi's
work and the fundamental work of Chern and Moser on the local
equivalence problem for real hypersurfaces\cite{CM74}. For background on real submanifolds in complex space, CR structures,
and their mappings, we refer to \cite{BER}. The projective-geometric viewpoint in CR geometry was subsequently
developed further by Burns and Shnider, who studied projective
connections naturally associated with CR structures in 
\cite{BurnsShnider80}.  A further fundamental development was Tanaka's formulation of the
geometry of Levi-nondegenerate real hypersurfaces in terms of graded
Lie algebras and canonical Cartan connections \cite{Tanaka76}, a
viewpoint that became one of the foundations of the modern theory of
CR and parabolic geometries. The correspondence between systems of partial differential equations,
their Lie symmetries, and CR geometry has been developed systematically
in modern form by Merker; see \cite{MerkerLie}.

\subsection{Real hypersurfaces and CR structures}

Let $
M\subset\mathbb{C}^{n+1} $
be a smooth real hypersurface. Write complex coordinates as
\[
(z,w)=(z^1,\ldots,z^n,w),
\]
where
\[
z=(z^1,\ldots,z^n)\in\mathbb{C}^n,
\qquad
w=u+iv\in\mathbb{C}.
\]
Locally, \(M\) may be defined by an equation
\[
\rho(z,w,\overline z,\overline w)=0,
\]
where \(\rho\) is real-valued and
\[
d\rho\neq0.
\]
The complex tangent space of \(M\) at \(p\) is
\[
H_pM
=
T_pM\cap J(T_pM),
\]
where \(J\) denotes multiplication by \(i\).
Equivalently, after complexification one considers
\[
T_p^{1,0}M
=
T_p^{1,0}\mathbb{C}^{n+1}
\cap
(T_pM\otimes\mathbb{C}).
\]
This is an \(n\)-dimensional complex vector space.

Thus a real hypersurface in \(\mathbb{C}^{n+1}\) naturally carries a
CR structure of CR dimension \(n\) and real dimension $
2n+1.
$
This dimension
should recall for us  the first jet space occurring in Hachtroudi's theory.

\subsection{The Levi form}

The first fundamental invariant of a CR hypersurface is the Levi form.
If
\[
L,L'\in T^{1,0}M,
\]
then, up to the usual choice of contact form, the Levi form is given by
\[
\mathcal{L}(L,\overline{L'})
=
-i\,d\theta(L,\overline{L'}),
\]
where \(\theta\) annihilates the complex tangent distribution.
In local defining coordinates one may equivalently write
\[
\mathcal{L}_{\alpha\overline{\beta}}
=
\frac{\partial^2\rho}
{\partial z^\alpha\partial\overline z^\beta}
\]
after restriction to complex tangent vectors.
The hypersurface is called \emph{Levi nondegenerate} if
\[
\det
\left(
\mathcal{L}_{\alpha\overline{\beta}}
\right)
\neq0.
\]

This condition is the CR analogue of nondegeneracy in several familiar
geometric settings. It is precisely the hypothesis under which the
Cartan--Chern--Moser theory has its cleanest form. The local equivalence problem for Levi-nondegenerate real
hypersurfaces in \(\mathbb{C}^{2}\) had already been studied
fundamentally by Cartan in his work on pseudo-conformal geometry
\cite{Cartan32I, Cartan32II}.  Cartan's approach provided a canonical
geometric structure whose curvature detects local equivalence with
the homogeneous model, and became one of the foundations of what is
now called CR geometry.

\subsection{The flat CR model}

The model for Levi-nondegenerate CR geometry is a nondegenerate
hyperquadric.
In the strongly pseudoconvex case one may take the unit sphere
\[
S^{2n+1}
=
\left\{
(z,w)\in\mathbb{C}^{n+1}:
|z|^2+|w|^2=1
\right\}.
\]
Removing one point and applying a Cayley transformation gives the
Heisenberg model
\[
\operatorname{Im}w
=
|z^1|^2+\cdots+|z^n|^2.
\]
More generally, if the Levi form has signature \((p,q)\), the model is
the hyperquadric
\[
\operatorname{Im}w
=
|z^1|^2+\cdots+|z^p|^2
-
|z^{p+1}|^2-\cdots-|z^n|^2.
\]
Now the CR equivalence problem asks:

\begin{quote}
When is a Levi-nondegenerate real hypersurface locally biholomorphically
equivalent to its model hyperquadric?
\end{quote}
This question is the CR counterpart of the question studied by
Hachtroudi which we recall again: 

\begin{quote}
When is a completely integrable second-order system locally equivalent
to the flat system \(y_{ij}=0\)?
\end{quote}
The surprising fact is that these are essentially the same question
after passing through Segre varieties.

\subsection{Real analyticity and complexification}

For this purpose one assumes that \(M\) is real analytic.
Locally one may solve its defining equation for \(w\) and write
\[
w
=
\Theta(z,\overline z,\overline w),
\]
where \(\Theta\) is holomorphic after the barred variables are regarded
as independent complex variables.
It is therefore useful to introduce independent complex parameters
\[
\zeta=(\zeta^1,\ldots,\zeta^n),
\qquad
\xi,
\]
in place of
\[
\overline z,\qquad\overline w.
\]
The complexified defining equation becomes
\[
w=\Theta(z,\zeta,\xi).
\]
For each fixed pair \((\zeta,\xi)\), this equation defines a complex
\(n\)-dimensional hypersurface
\[
S_{\zeta,\xi}
=
\left\{
(z,w):
w=\Theta(z,\zeta,\xi)
\right\}.
\]
These hypersurfaces are the \emph{Segre varieties} of \(M\).
Thus a real analytic CR hypersurface naturally determines an
\((n+1)\)-parameter family of complex hypersurfaces.
This should now be compared  with the general solution
\[
y=Q(x,a,b)
\]
of the completely integrable differential systems  treated by Hachtroudi and considered earlier.

\subsection{From Segre varieties to a differential equation}

The connection with Hachtroudi's theory now becomes explicit.
Start with
\[
w=\Theta(z,\zeta,\xi).
\]
Differentiate with respect to \(z^i\):
\[
w_{z^i}
=
\Theta_{z^i}(z,\zeta,\xi),
\qquad
1\leq i\leq n.
\]
Together with the original equation, we obtain \(n+1\) equations
\[
w=\Theta(z,\zeta,\xi),
\]
\[
w_{z^1}=\Theta_{z^1}(z,\zeta,\xi),
\quad\ldots,\quad
w_{z^n}=\Theta_{z^n}(z,\zeta,\xi).
\]
Under the Levi-nondegeneracy assumption these equations can locally be
solved for the \(n+1\) parameters $
(\zeta,\xi) $
in terms of $
(z,w,w_z). $
Thus
\[
\zeta
=
\zeta(z,w,w_z),
\qquad
\xi
=
\xi(z,w,w_z).
\]
If we differentiate  one more time 
\[
w_{z^iz^j}
=
\Theta_{z^iz^j}(z,\zeta,\xi).
\]
and substituting the expressions for \(\zeta\) and \(\xi\)  we obtain
\[
w_{z^iz^j}
=
\Phi_{ij}(z,w,w_z),
\qquad
1\leq i,j\leq n.
\]
We have therefore obtained a completely integrable system of the form
\[
w_{ij}
=
\Phi_{ij}(z,w,w_z),
\]
which is precisely the class of systems studied by Hachtroudi.
The general solution of this system is exactly the Segre family
\[
w=\Theta(z,\zeta,\xi).
\]
Thus we have the fundamental passage

\[
\text{Real analytic CR hypersurface}
\longrightarrow
\text{Segre family}\]
\[ \longrightarrow
\text{Completely integrable second-order PDE}.
\]
This is the central bridge between Hachtroudi's geometry and CR
geometry.

\subsection{Why Levi nondegeneracy appears}

It is worth emphasizing the role of the Levi form in this construction.
To eliminate the parameters \((\zeta,\xi)\), one must invert the map
\[
(\zeta,\xi)
\longmapsto
\left(
\Theta,
\Theta_{z^1},
\ldots,
\Theta_{z^n}
\right).
\]
Its Jacobian determinant is
\[
\Delta
=
\det
\begin{pmatrix}
\Theta_{\zeta^1} & \cdots &
\Theta_{\zeta^n} & \Theta_{\xi}
\\
\Theta_{z^1\zeta^1} & \cdots &
\Theta_{z^1\zeta^n} & \Theta_{z^1\xi}
\\
\vdots && \vdots & \vdots
\\
\Theta_{z^n\zeta^1} & \cdots &
\Theta_{z^n\zeta^n} & \Theta_{z^n\xi}
\end{pmatrix}.
\]
The nonvanishing condition
\[
\Delta\neq0
\]
is equivalent, in suitable coordinates, to Levi nondegeneracy.

Thus the very hypothesis which is fundamental in CR geometry is also
exactly what permits the elimination procedure leading to
Hachtroudi's differential equation.
This is not an accidental analogy: the two theories are describing
the same underlying geometry in two different languages.

\subsection{The sphere and the flat differential equation}

Consider the Heisenberg hypersurface
\[
\operatorname{Im}w
=
|z|^2.
\]
After complexification, its Segre varieties are affine complex
hyperplanes.
In suitable coordinates they have the form
\[
w
=
b+\sum_{i=1}^{n}a_i z^i,
\]
where \((a_1,\ldots,a_n,b)\) are parameters.
Consequently,
\[
w_{z^iz^j}=0.
\]
We see that  the CR-flat model corresponds exactly to Hachtroudi's flat model:
\[
\text{hyperquadric}
\quad\longleftrightarrow\quad
w_{z^iz^j}=0.
\]
It follows that a biholomorphic transformation sending a
Levi-nondegenerate hypersurface to a hyperquadric simultaneously
straightens its Segre varieties into affine hyperplanes.

From the PDE point of view, this means that the associated system is
transformed into $
w_{z^iz^j}=0.
$

Hence Hachtroudi's flatness problem and the CR sphericity problem are
two manifestations of the same equivalence problem.

\subsection{The Hachtroudi curvature}

Assume throughout this subsection that $n\geq 2$. For a completely integrable system
\[
y_{ij}=F_{ij}(x,y,p),
\]
Hachtroudi's curvature contains a particularly important trace-free
component involving the second derivatives of \(F_{ij}\) with respect
to the first-derivative variables \(p_k\).

A convenient way of writing its vanishing condition is
\[
\begin{split}
0={}&
\frac{\partial^2F_{ij}}
{\partial p_k\partial p_\ell}
\\
&-\frac{1}{n+2}
\sum_{r=1}^n
\left(
\delta_i^k
\frac{\partial^2F_{rj}}
{\partial p_r\partial p_\ell}
+
\delta_i^\ell
\frac{\partial^2F_{rj}}
{\partial p_k\partial p_r}
\right.
\\
&\hspace{35mm}\left.
+
\delta_j^k
\frac{\partial^2F_{ir}}
{\partial p_r\partial p_\ell}
+
\delta_j^\ell
\frac{\partial^2F_{ir}}
{\partial p_k\partial p_r}
\right)
\\
&+
\frac{
\delta_i^k\delta_j^\ell+
\delta_j^k\delta_i^\ell
}
{(n+1)(n+2)}
\sum_{r,s=1}^n
\frac{\partial^2F_{rs}}
{\partial p_r\partial p_s}.
\end{split}
\]
The important point is not the appearance of this formula itself, but
its geometric meaning.
It is the trace-free curvature of the canonical projective geometry
associated with the system.

For \(n=1\), the displayed trace-free tensor vanishes identically for every
scalar second-order equation. Put
\[
A=\frac{\partial^2 F_{11}}{\partial p_1^2}.
\]
All indices and all Kronecker deltas are equal to \(1\), so the displayed
expression reduces to
\[
A-\frac{1}{3}(4A)+\frac{2}{(2)(3)}A
=
A-\frac{4}{3}A+\frac{1}{3}A
=
0.
\]
Thus this tensor cannot detect flatness in dimension one.

Let us give a concrete  counterexample.
Consider the scalar ordinary differential equation
\[
y''=(y')^4.
\]
For \(n=1\), complete integrability is automatic: the associated distribution
has rank one. The Hachtroudi  tensor vanishes identically by the
preceding calculation. Nvertheless, the equation is not locally point-equivalent to
\[
Y''=0.
\]
Indeed, let
\[
X=X(x,y),
\qquad
Y=Y(x,y)
\]
be a local point transformation, write
\[
p=y',
\]
and suppose that the image equation is \(Y''=0\). The transformed first
derivative is
\[
P=\frac{dY}{dX}
=
\frac{Y_x+pY_y}{X_x+pX_y}.
\]
Writing
\[
D=\partial_x+p\,\partial_y+F\,\partial_p,
\]
the condition \(Y''=0\) is equivalent to
\[
D(P)=0.
\]
Solving this equation for \(F\) gives
\[
F=
-\frac{
\bigl(Y_{xx}+2pY_{xy}+p^2Y_{yy}\bigr)
\bigl(X_x+pX_y\bigr)
-
\bigl(Y_x+pY_y\bigr)
\bigl(X_{xx}+2pX_{xy}+p^2X_{yy}\bigr)
}{
Y_yX_x-Y_xX_y
}.
\]
The denominator is the nonzero Jacobian of the point transformation, and the
numerator is a polynomial of degree at most three in \(p\). Hence every point
transform of \(Y''=0\) has a right-hand side that is cubic, or of lower
degree, in \(p\). The equation
\[
y''=p^4
\]
therefore cannot be locally point-equivalent to \(Y''=0\).

Consequently, in dimension \(n=1\), the displayed trace-free tensor vanishes
identically even for nonflat equations. The flatness problem must instead be
detected by the higher-order Cartan--Tresse invariants.

For the differential system obtained from the Segre family of a
Levi-nondegenerate CR hypersurface, this tensor becomes the
Hachtroudi--Chern--Moser curvature,  and  its vanishing means that the Segre family can locally be straightened
into affine hyperplanes and hence that the original CR hypersurface is
locally equivalent to the appropriate hyperquadric.
Schematically,
\[
\text{Hachtroudi curvature}=0
\Longleftrightarrow
\text{PDE is projectively flat}\]
\[ \Longleftrightarrow
\text{Segre varieties are affine hyperplanes}\]
\[ \Longleftrightarrow 
\text{CR hypersurface is locally pseudospherical}.
\]
For CR dimension \(n\geq2\), this curvature is closely related to what
is now usually called the Chern--Moser or Chern--Moser--Weyl tensor.

\subsection{Chern and Moser}

In 1974 Shiing-Shen Chern and J\"urgen Moser published their
fundamental paper
\emph{Real hypersurfaces in complex manifolds}\cite{CM74}.
Their problem was essentially the several-complex-variable version of
Poincar\'e's equivalence problem:

\begin{quote}
When are two Levi-nondegenerate real hypersurfaces locally equivalent
under a biholomorphic transformation?
\end{quote}

The paper solved this problem by combining Cartan's geometric ideas
with a powerful normal-form construction.
Near a point of a Levi-nondegenerate hypersurface, suitable
holomorphic coordinates put the defining equation into the form
\[
v
=
\langle z,\overline z\rangle
+
F(z,\overline z,u),
\]
where
\[
\langle z,\overline z\rangle
=
\sum_{\alpha=1}^{n}
\epsilon_\alpha
z^\alpha\overline z^\alpha,
\qquad
\epsilon_\alpha=\pm1.
\]
The higher-order term \(F\) can be subjected to normalization
conditions which remove precisely the coordinate freedom coming from
local biholomorphic transformations.
In the strongly pseudoconvex case the model term is simply
\[
v=|z|^2.
\]
The first essential trace-free term remaining after normalization
produces the Chern--Moser curvature tensor.
Thus the philosophy is again Cartan's:
\[
\text{geometry}
\longrightarrow
\text{normalization}
\longrightarrow
\text{canonical curvature}
\longrightarrow
\text{equivalence criterion}.
\]

\subsection{Hachtroudi and the Chern--Moser connection}

The historical relationship between these constructions deserves
special emphasis.
Shortly after the Chern--Moser paper, Chern published a note entitled
\emph{On the projective structure of a real hypersurface in
\(\mathbb{C}^{n+1}\)}\cite{Chern75}.
Chern begins from Segre's idea of replacing the conjugate variables in
the defining equation of a real analytic hypersurface by independent
complex parameters. This produces the family of Segre varieties
described above.
He then observes that this family naturally defines a projective
connection in the space of hyperplane elements.

At this point Chern explicitly invokes Hachtroudi's work. He explains
that, apart from notation, the projective connection obtained in this
way is essentially the connection constructed by Hachtroudi.
Even more importantly, Chern shows that this projective connection is
closely related to, and in fact underlies, the connection constructed
in the joint work of Chern and Moser.
Thus the historical and mathematical chain is remarkably direct:
\[
\text{Cartan's projective geometry}
\longrightarrow
\text{Hachtroudi, 1937}
\longrightarrow \]
\[ \text{Segre--Hachtroudi projective connection}
\longrightarrow
\text{Chern--Moser CR geometry}.
\]

This provides perhaps the clearest evidence of the lasting
mathematical importance of Hachtroudi's thesis.

\subsection{Two approaches to the same curvature}

There are consequently two different ways of reaching essentially the
same local CR invariant.
The first is the Chern--Moser route:
\[
M
\longrightarrow
\text{normal form}
\longrightarrow
\text{Chern--Moser tensor}.
\]
The second is the Hachtroudi route:
\[
M
\longrightarrow
\text{Segre varieties}
\longrightarrow
w_{ij}=\Phi_{ij}(z,w,w_z)
\longrightarrow
\text{Hachtroudi curvature}.
\]
The first works directly with the real hypersurface and its
biholomorphic normalization.
The second passes through the holomorphic geometry of the Segre
family and interprets the problem as the projective equivalence of a
system of differential equations.

Modern treatments show that these two curvatures encode the same
flatness obstruction.
This is a particularly beautiful example of how a construction
originating in the geometry of differential equations in the 1930s
reappeared, almost forty years later, at the heart of several complex
variables and CR geometry.

\subsection{The special case of CR dimension one}

There is one important qualification.
When $
n=1, $
so that \(M\) is a three-dimensional real hypersurface in
\(\mathbb{C}^2\), the lowest-order Chern--Moser tensor which appears in
higher CR dimension vanishes identically for dimensional reasons.

The fundamental curvature invariant is then of higher order. This is
the geometry already studied by Cartan in his 1932-1933 work on
three-dimensional CR structures\cite{Cartan32I, Cartan32II}.
Thus one should distinguish $
n=1
$
from
$
n\geq2.
$
For \(n\geq2\), the trace-free Hachtroudi curvature and the
Chern--Moser curvature give the natural lowest-order obstruction to
sphericity.
For \(n=1\), the corresponding obstruction is Cartan's higher-order
CR curvature.
This dimensional distinction is important both historically and
mathematically.

\subsection{ A modern interpretation}

From the modern viewpoint the relation between the two theories can
be summarized in terms of Cartan geometries.
A Levi-nondegenerate CR hypersurface carries a canonical parabolic
geometry whose flat model is the appropriate hyperquadric.
The Segre family provides a holomorphic projective realization of
this geometry.

Hachtroudi's normal projective connection is therefore not merely an
analogy with the Chern--Moser connection. It is one of the geometric
mechanisms through which the CR curvature may be recovered.

The passage
\[
\text{CR geometry}
\longleftrightarrow
\text{Segre geometry}
\longleftrightarrow
\text{PDE geometry}
\]
is consequently one of the most important modern interpretations of
Hachtroudi's work.
It reveals that his 1937 thesis anticipated structures which later
became central in Cartan geometry, CR geometry, and the modern
geometric theory of differential equations.

\medskip

\section{Chern's Projective Structure and the Legacy of Hachtroudi}

\begin{quote}
{\it The purpose of this note is to carry out Segre's idea for general
\(n\) and relate it to the invariants given in \cite{CM74}. Equation
(2) defines an \((n+1)\)-parameter family of hypersurfaces when
\(M\) is non-degenerate. Generalizing the work of Arthur Tresse,
M.~Hachtroudi showed that a projective connection can be defined
intrinsically in the space of hyperplane elements of
\(\mathbf{C}^{n+1}\) \cite{Hach37}. The definition is a
generalization, by no means obvious, of the construction of classical
projective geometry from the data of its hyperplanes; $\dots$ We will
show that the definition of Hachtroudi's connection is closely related
to that of the connection in \cite{CM74}.}
\end{quote}

\hfill---S.~S.~Chern \cite[p.~74]{Chern75}

\medskip

\begin{figure}[H]
   \centering
   \includegraphics[width=0.6\textwidth]{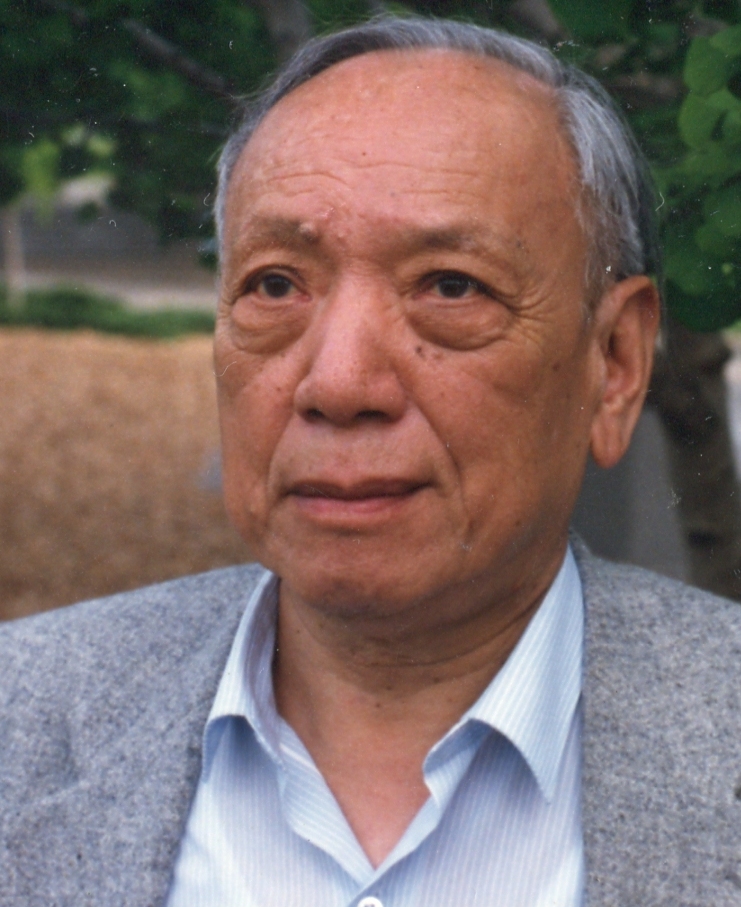} 
   \caption{Shiing-Shen Chern (26 October 1911--3 December 2004).
He received his doctorate from the University of Hamburg in 1936 under
Wilhelm Blaschke and spent the academic year 1936--1937 in Paris working
with \'{E}lie Cartan. During this period he overlapped in Cartan's circle
with Hachtroudi, who was completing his own doctorate in 1937; it was
probably then that their acquaintance began. Hachtroudi later listed Chern
among the mathematicians with whom he maintained correspondence and
scientific relations.}
   \end{figure}

As we indicated earlier, the relation between Hachtroudi's work and CR geometry became
particularly transparent in S.~S.~Chern's 1975 paper
\emph{On the projective structure of a real hypersurface in
\(\mathbb{C}^{n+1}\)} \cite{Chern75}. The paper is important for our
purposes for two reasons. Mathematically, Chern showed how the Segre
family of a Levi-nondegenerate real hypersurface naturally gives rise
to a projective connection. Historically, he explicitly identified
this construction with the projective connection introduced by
Hachtroudi almost forty years earlier and explained its close relation
with the connection arising in the work of Chern and Moser.

Chern's observation provides one of the clearest indications of the
lasting significance of Hachtroudi's thesis. A construction originally
developed in connection with the equivalence problem for systems of
second-order partial differential equations reappears naturally in the
local geometry of real hypersurfaces in complex space.

\subsection{Segre families and differential equations}

Let
\[
M\subset\mathbb{C}^{n+1}
\]
be a real-analytic hypersurface, locally defined by
\[
r(z,\overline z)=0.
\]
Segre's idea is to replace the conjugate variables \(\overline z\) by
independent complex parameters \(a\) and consider
\[
r(z,a)=0.
\]
For each fixed value of \(a\), this equation defines a complex
hypersurface. In this way \(M\) determines an
\((n+1)\)-parameter family of complex hypersurfaces, its
\emph{Segre family}.

Locally, write the variables as
\[
(z^1,\ldots,z^n,w)\in\mathbb{C}^{n+1}.
\]
Under the Levi-nondegeneracy assumption, the Segre family can be used
to produce a completely integrable system of second-order differential
equations
\[
w_{\alpha\beta}
=
F_{\alpha\beta}(z,w,w_\gamma),
\qquad
1\leq\alpha,\beta\leq n.
\]
Its solutions are precisely the members of the Segre family. Thus a
problem in CR geometry is converted into a problem concerning the
geometry of a system of differential equations.

This is exactly the setting in which Hachtroudi's construction becomes
relevant. His 1937 thesis associates a canonical projective connection
with a completely integrable system of this type. The passage
\[
M
\longrightarrow
\text{Segre family}
\longrightarrow
\text{second-order PDE}
\]
therefore provides a direct bridge between CR geometry and
Hachtroudi's geometry of differential equations. For the role of Segre
families in this context, see for example \cite{Faran80,Sukhov01}.

\subsection{Hyperplane elements and contact geometry}

There is a second, geometric way of seeing the same construction.
A hyperplane element in \(\mathbb{C}^{n+1}\) consists of a point
together with a complex hyperplane through that point. If a
hypersurface is locally represented as
\[
w=w(z^1,\ldots,z^n),
\]
its tangent hyperplane is determined by
\[
p_\alpha=\frac{\partial w}{\partial z^\alpha}.
\]
Consequently the space of hyperplane elements has local coordinates
\[
(z^1,\ldots,z^n,w,p_1,\ldots,p_n).
\]
This is precisely the local coordinate description of the first jet
space
\[
J^1(\mathbb{C}^n,\mathbb{C}).
\]
Thus the two notions that occur naturally in the two theories are
really the same geometric object:
\[
\text{hyperplane element}
\quad\longleftrightarrow\quad
\text{first jet}.
\]

On this space there is the canonical contact form
\[
\theta
=
dw-\sum_{\alpha=1}^{n}p_\alpha\,dz^\alpha.
\]
The equation \(\theta=0\) expresses the condition that a hyperplane
element be tangent to a hypersurface. Its kernel
\[
H=\ker\theta
\]
is a \(2n\)-dimensional contact distribution.

The Segre family determines one distinguished family of integral
submanifolds of this contact structure, while the fibers obtained by
fixing the base point determine another. This is the same double
incidence geometry that appeared earlier in our discussion of
Hachtroudi's differential system. Chern's space of hyperplane elements
and Hachtroudi's first jet space are therefore not merely analogous:
they are two descriptions of the same underlying contact geometry.

\subsection{The projective connection}

Chern proceeds by choosing a coframe adapted to this contact geometry.
The basic forms may be written schematically as
\[
\theta,\qquad
\theta^\alpha,\qquad
\theta_\alpha,
\]
with the contact structure encoded by a structure equation of the form
\[
d\theta
=
i\,\theta^\alpha\wedge\theta_\alpha
+
\theta\wedge\varphi.
\]
Further connection forms are introduced by differentiating the
structure equations and imposing normalization conditions. The
resulting collection of forms constitutes a canonical projective
connection.

The complexified Segre geometry has symmetry algebra
\[
\mathfrak{sl}(n+2,\mathbb{C}).
\]
This is the same projective symmetry algebra that
appears in the flat model of Hachtroudi's geometry. The occurrence of
this Lie algebra in both constructions is  intrinsic rather
than accidental. We note that the corresponding real CR hyperquadric of Levi signature \((p,q)\) has
symmetry algebra
\[
\mathfrak{su}(p+1,q+1).
\]

At this point Chern makes the crucial identification. Apart from
differences in notation, the projective connection obtained from the
Segre family is essentially the connection constructed by Hachtroudi.
The relation may be summarized as
\[
\text{Segre family}
\longrightarrow
\text{Space of hyperplane elements}\]
\[ \longrightarrow
\text{Hachtroudi projective connection}.
\]

Chern then compares these structure equations with those used in his
joint work with Moser. He shows that they can be normalized so as to
agree with the corresponding Chern--Moser structure equations
\cite{CM74,Chern75}. Thus Hachtroudi's projective construction is
closely related to the canonical geometry underlying the
Chern--Moser theory of Levi-nondegenerate real hypersurfaces.

This is perhaps the central mathematical point of Chern's 1975 paper
for the present article. It establishes a direct connection between
Hachtroudi's work on differential equations and one of the fundamental
constructions of modern CR geometry.

\subsection{Curvature and the flat model}

As in Cartan's general approach to equivalence problems, the decisive
invariant is curvature. If \(\Pi\) denotes the matrix-valued
projective connection form, its curvature is
\[
K=d\Pi+\Pi\wedge\Pi.
\]
For the homogeneous projective model,
\[
K=0.
\]
The curvature therefore measures the obstruction to local equivalence
with the flat projective incidence geometry.

In terms of Hachtroudi's differential equations, the flat model is
\[
w_{\alpha\beta}=0.
\]
Its solutions are affine hyperplanes, and the corresponding incidence
geometry is the homogeneous point--hyperplane model. Consequently,
\[
K=0 \quad \text{locally}
\quad\Longleftrightarrow\quad
\text{local equivalence with }
w_{\alpha\beta}=0.
\]

When the differential system arises from the Segre family of a
Levi-nondegenerate real hypersurface, the same condition has a CR
interpretation: the hypersurface is locally CR equivalent to a
nondegenerate hyperquadric. Thus the various descriptions of flatness
fit into a single picture:
\[
\begin{array}{c}
\text{vanishing projective curvature}
\\[1mm]
\Updownarrow
\\[1mm]
\text{flat Hachtroudi geometry}
\\[1mm]
\Updownarrow
\\[1mm]
w_{\alpha\beta}=0
\\[1mm]
\Updownarrow
\\[1mm]
\text{flat Segre geometry}
\\[1mm]
\Updownarrow
\\[1mm]
\text{CR hyperquadric}.
\end{array}
\]

For CR dimension \(n\geq2\), the fundamental CR curvature is represented
by the trace-free Chern--Moser tensor
\[
S_{\alpha\overline{\beta}\gamma\overline{\delta}}.
\]
Its vanishing characterizes local CR equivalence with the hyperquadric.
From the Segre--Hachtroudi viewpoint, this curvature is therefore not
an unrelated invariant: it expresses, in CR language, the obstruction
to flatness of the projective geometry associated with the Segre
family.

There is one important low-dimensional qualification. When \(n=1\),
so that \(M\) has real dimension three, the Chern--Moser tensor
vanishes identically for algebraic reasons. The first nontrivial
curvature occurs at higher order, in accordance with Cartan's
classical treatment of three-dimensional CR geometry. Chern's 1975
analysis also reflects this exceptional behavior \cite{Chern75}.  The Cartan equivalence problem for Levi-nondegenerate real
hypersurfaces has been revisited in modern form by Jo\"el Merker and
Masoud Sabzevari \cite{MerkerSabzevari12,MerkerSabzevari14}, including
an explicit construction of Cartan's connection for
Levi-nondegenerate real hypersurfaces in \(\mathbb{C}^{2}\).

\subsection{From Hachtroudi to modern geometry}

In this subsection and the next  we  summarize Hachtroudi's contribution to the geometry of differential equations as we discussed in this paper.  We showed that Hachtroudi's contribution belongs to a long development in the
geometry of differential equations. Its origins lie in the
nineteenth-century work of Lie and Tresse on transformation groups and
differential invariants \cite{LieEngel,Tresse}. Cartan profoundly
changed the nature of the equivalence problem by replacing the search
for individual invariants with the construction of canonical
geometric structures and their curvature.

Hachtroudi applied this philosophy systematically to completely
integrable systems of second-order partial differential equations.
His central achievement was to associate with such a system a
canonical projective connection whose curvature measures the
obstruction to equivalence with the flat system
\[
w_{\alpha\beta}=0.
\]
In this sense he transformed an equivalence problem for nonlinear
partial differential equations into a curvature problem in projective
differential geometry.

Almost four decades later, the same geometry reappeared from a
different direction. Chern and Moser developed their theory of
Levi-nondegenerate real hypersurfaces in 1974 \cite{CM74}, and in the
following year Chern explicitly returned to the projective viewpoint
and identified the connection associated with the Segre family with
Hachtroudi's construction \cite{Chern75}. One may therefore summarize
the historical development as
\[
\text{Lie--Tresse}
\longrightarrow
\text{Cartan}
\longrightarrow
\text{Hachtroudi}
\longrightarrow\]
\[ \text{Chern--Moser}
\longrightarrow
\text{Modern CR geometry}.
\]

There is also a natural modern interpretation. The flat
point--hyperplane incidence geometry is a homogeneous geometry for
\(SL(n+2,\mathbb{C})\), and its curved analogues belong naturally to
the general framework of Cartan and parabolic geometries
\cite{CapSlovak}. In this language a curved geometry is equipped with
a Cartan connection
\[
\omega\in\Omega^1(\mathcal{G},\mathfrak{g}),
\]
whose curvature
\[
\Omega
=
d\omega+\frac12[\omega,\omega]
\]
measures the obstruction to local equivalence with the homogeneous
model. The terminology and general theory of parabolic geometry were
developed much later, but the essential mechanism is already visible
in Hachtroudi's thesis.

This modern viewpoint also helps explain the continued appearance of
Hachtroudi's name in the literature. Jet spaces, contact structures,
integrable distributions, canonical Cartan connections, and curvature
as an obstruction to flatness are now standard ingredients in the
geometric study of differential equations. In Hachtroudi's work these
ideas already occur together in a remarkably concrete form.

\subsection{Concluding perspective}

Chern's 1975 paper brings the story full circle. Starting with a
Levi-nondegenerate real hypersurface, one passes through its Segre
family to a system of differential equations and hence to
Hachtroudi's projective connection:
\[
M
\longrightarrow
\text{Segre family}
\longrightarrow
\text{Second-order PDE}\]
\[ \longrightarrow
\text{Hachtroudi connection}
\longrightarrow
\text{Curvature}.
\]
The Chern--Moser approach begins with the same hypersurface and arrives,
through a different route, at its canonical CR curvature. Chern's
observation shows that these two approaches are intimately connected.

Hachtroudi's thesis should therefore not be viewed merely as an
isolated contribution to the theory of differential equations.
Written by a young Iranian mathematician working under \'Elie Cartan
in Paris in the 1930s, it belongs naturally to a line of ideas that
runs from Lie and Tresse through Cartan to Chern and Moser, and
continues today in Cartan, CR, and parabolic geometry.

Perhaps the most enduring lesson of Hachtroudi's construction is a
simple one: nonlinear differential equations may themselves carry
canonical geometric structures, and their curvature reveals the
obstruction to reducing them to a flat model. This principle, so
familiar in modern differential geometry, is already clearly present
in Hachtroudi's work.  The geometric meaning of the vanishing of Hachtroudi's curvature has
also been revisited explicitly in recent work of Merker, relating the
vanishing curvature condition to local equivalence with the flat CR
model \cite{MerkerHach}.

\appendix

\section{Hachtroudi in His Own Words}

The following autobiographical sketch was written by Mohsen
Hachtroudi in Persian for the mathematics magazine
\emph{Yek\=an} in 1964 \cite{HachAuto}. 
It provides  a rare first-hand account of
Hachtroudi's education, his teachers, his mathematical work, and his
scientific career. Equally revealing is the breadth of his
intellectual interests, extending well beyond mathematics to
education, literature, and culture; qualities that were repeatedly
emphasized by his students, colleagues, and contemporaries.

\medskip 

\begin{quote}
\noindent
{\it I was born in Tabriz on 22 Dey 1286 (12 January 1908)\footnote{For the birth date conversion see the
 footnote in page 5.}.  I completed my
primary education at the Aghdasieh and Sirus schools, and my secondary
education at D\=ar al-Fon\=un. In 1304 (1925) I graduated from
D\=ar al-Fon\=un and then studied medicine for several years.

My first trip to Europe was followed by a return to Iran. Upon my
return, I enrolled in the Higher Teachers' College
(\emph{D\=ar al-Mo`allem\={\i}n-e Markaz\={\i}}), which had recently
been established, where I chose mathematics as my field of study. I
graduated as a member of its second graduating class.

I then made a second trip to France, where I studied at the Faculty of
Sciences of the University of Paris. I obtained the degree of
\emph{Licence} and, in 1937, received the degree of
\emph{Doctorat d'\'Etat} from the Sorbonne.

In 1315 (1937) I returned to Iran and was appointed Associate Professor
in the Faculty of Sciences and at the Higher Teachers' College. In
1320 (1942) I was promoted to the rank of Professor. Thereafter, in
addition to my university professorship, I held the following
positions. In 1321 (1943) I served as Director of Education for
Tehran; in 1330 (1952) I became President of the University of Tabriz;
and in 1336 (1958) I served for one term as Dean of the Faculty of
Sciences of the University of Tehran.

In 1323 (1945) I married. I have three children: one son and two
daughters.

\noindent
Among my teachers, first and foremost was my late brother,
Mohammad Zia Hachtroudi, to whom I owe the very foundations of my
education. The late Gholamhossein Rahnam\=a was another of my
teachers; indeed, many of today's professors and school teachers
have, directly or indirectly, been trained by him.

Among those teachers who are still living, all have contributed to my
education, but I must especially mention Professor
Abdolazim Gharib of the University of Tehran. Whatever knowledge and
literary cultivation I possess, I owe largely to him. Even more
important was the moral and spiritual education that he gave to many
of us students. Spiritually and ethically I remain indebted to all my
teachers and mentors.

Mr.  Kayhan, besides being one of my teachers, also watched over
my education for many years during his service as Vice-President of
the University.

Outside the teaching profession, many of my successes (which in truth
have been modest and insignificant) are due to the encouragement and
support of my mentor, Dr.\ Siy\=asi. Everyone knows how much he
contributed to the establishment of the University of Tehran and to
the respect and dignity accorded to university professors. For this,
the academic community owes him a lasting debt. Even today he remains
a guide and source of inspiration for those, like myself, who seek
knowledge.

Among my university professors, first and foremost was the late
Professor \'Elie Cartan, the founder of modern geometry, to whom
almost every branch of contemporary mathematics is indebted. He was
the supervisor of my doctoral dissertation and, more generally, my
scientific mentor.

The present professors of the University of Paris, most of whom were
students of Professor Cartan, were my contemporaries and also guided
my work. Among them were Professor Ehresmann at the Sorbonne,
Professor Lichnerowicz at the Coll\`ege de France, the late Professor
Weyl at Princeton University in the United States (who was my teacher
and adviser during a period of study that I spent at Princeton),
Professor Schouten in Delft, Holland (now Director of the Amsterdam
Mathematical Centre), Professor Vinogradov and the late Professor
Finikov of Moscow University, together with several others.

The subject of my doctoral dissertation concerned projective spaces
of elements (a point together with a line or a plane) endowed with
normal connections, and, as stated earlier, my supervisor was the late
Professor \'Elie Cartan.

Among my writings and articles in mathematics are several papers on
the infinitesimal geometry of general spaces, especially non-holonomic
spaces; several papers on the geometry of normal spaces and figures
with internal similarity; several papers on analytical mechanics and
dynamical trajectories; several papers on differential equations; a
generalization of Euler's theorems on continued fractions; the
application of continued fractions to the solution of differential
equations, and in particular to the determination of the solvable
cases of the Riccati equation; the law of dual algebra in the mechanics
of higher-dimensional spaces; Weyl spaces with orthogonal connection;
Schouten spaces with constant invariants; and several other papers.
Some of these articles have been published, in the form of collections
in the French language, by the University of Tehran.

Up to the present time I have participated in four International
Congresses of Mathematicians: Harvard and Cambridge in the United
States, Amsterdam in Holland, and Edinburgh in England, as well as
Nice (the Congress of Mathematicians of the Latin Languages).

In addition to participating in congresses, I have attended and
lectured at scientific gatherings in various countries by invitation.
Among these were invitations from the Academy of Sciences of the
Soviet Union, the Academy of Sciences of Bucharest, Moscow University,
a second invitation from the Academy of Sciences of the Soviet Union,
the Rehovot Institute in Israel, where I was introduced to the
professors of that institute, the Scientific Congress of Pakistan,
and the Faculty of Sciences of Paris, which invited me to lecture.

Among the scholars with whom I have maintained correspondence and
scientific relations, in addition to the professors already mentioned,
are Professor Struik of M.I.T. in the United States; Professor
Zariski; Professor Chern of Chicago; Professor Albert in the United
States; Professor Ishlinskii of Moscow; Alexandrov of Leningrad;
Luzhin of Kiev; Javad Maghsudov of Baku; Segre and Bompiani of Italy;
several professors in London and Manchester; Professors Haimovici
(two brothers) of Bucharest; Moisil and Ronceanu of Bucharest;
Professor Stoilow, President of the Mathematical Academy of
Bucharest; and the Secretary-General of the Bucharest Academy.}
\end{quote}


\begin{thebibliography}{99}

\bibitem{AkivisRosenfeld}
M. A. Akivis and B. A. Rosenfeld,
\emph{\'Elie Cartan (1869--1951)},
Translated from the Russian by V. V. Goldberg,
Translations of Mathematical Monographs,
Vol.~123,
American Mathematical Society,
Providence, RI,
1993.

\bibitem{BER}
M. S. Baouendi, P. Ebenfelt and L. P. Rothschild,
\emph{Real Submanifolds in Complex Space and Their Mappings},
Princeton Mathematical Series, vol.~47,
Princeton University Press,
Princeton, NJ, 1999.

\bibitem{Bieche}
C. Bi\`eche,
Le probl\`eme d'\'equivalence locale pour un syst\`eme scalaire
complet d'\'equations aux d\'eriv\'ees partielles d'ordre deux
\`a $n$ variables ind\'ependantes,
\emph{Ann. Fac. Sci. Toulouse Math.}
{\bf 16} (2007), 1--36.


\bibitem{BCGGG}
R. L. Bryant, S. S. Chern, R. B. Gardner,
H. L. Goldschmidt and P. A. Griffiths,
\emph{Exterior Differential Systems},
Mathematical Sciences Research Institute Publications,
vol.~18,
Springer,
New York, 1991.


\bibitem{BurnsShnider80}
D. Burns and S. Shnider,
Projective connections in CR geometry,
\emph{Manuscripta Math.}
{\bf 33} (1980), 1--26.

\bibitem{Cartan24}
\'E. Cartan,
Sur les vari\'et\'es \`a connexion projective,
\emph{Bull. Soc. Math. France}
{\bf 52} (1924), 205--241.


\bibitem{Cartan32I}
\'E. Cartan,
Sur la g\'eom\'etrie pseudo-conforme des hypersurfaces de
l'espace de deux variables complexes. I,
\emph{Ann. Mat. Pura Appl.}
{\bf 11} (1932), 17--90.

\bibitem{Cartan32II}
\'E. Cartan,
Sur la g\'eom\'etrie pseudo-conforme des hypersurfaces de
l'espace de deux variables complexes. II,
\emph{Ann. Scuola Norm. Sup. Pisa}
{\bf 1} (1932), 333--354.



\bibitem{Cartan35}
\'E. Cartan,
\emph{La m\'ethode du rep\`ere mobile,
la th\'eorie des groupes continus et les espaces g\'en\'eralis\'es},
Actualit\'es Scientifiques et Industrielles,
Hermann,
Paris, 1935.

\bibitem{CartanNotice}
\'E. Cartan,
\emph{Notice sur les travaux scientifiques;
suivi de Le parall\'elisme absolu et la th\'eorie unitaire du champ},
Collection \emph{Discours de la m\'ethode},
No.~6,
Gauthier--Villars,
Paris, 1974.



\bibitem{CapSlovak}
A. \v{C}ap and J. Slov\'ak,
\emph{Parabolic Geometries I:
Background and General Theory},
Mathematical Surveys and Monographs,
vol.~154,
American Mathematical Society,
Providence, RI, 2009.

\bibitem{ChernChevalley}
S. S. Chern and C. Chevalley,
\'Elie Cartan and his mathematical work,
\emph{Bull. Amer. Math. Soc.}
{\bf 58} (1952), 217--250.




\bibitem{Chern75}
S. S. Chern,
On the projective structure of a real hypersurface in
$\mathbb{C}^{n+1}$,
\emph{Math. Scand.}
{\bf 36} (1975), 74--82.

\bibitem{CM74}
S. S. Chern and J. K. Moser,
Real hypersurfaces in complex manifolds,
\emph{Acta Math.}
{\bf 133} (1974), 219--271.


\bibitem{ChernSelectedI}
S.-S.~Chern,
\textit{Selected Papers, Vol.~I},
Springer-Verlag, New York, 1978.



\bibitem{Faran80}
J. J. Faran,
Segre families and real hypersurfaces,
\emph{Invent. Math.}
{\bf 60} (1980), 135--172.




\bibitem{Fels95}
M. E. Fels,
The equivalence problem for systems of second-order ordinary
differential equations,
\emph{Proc. London Math. Soc.}
{\bf 71} (1995), 221--240.




\bibitem{Grossman00}
D. A. Grossman,
Torsion-free path geometries and integrable second-order ODE systems,
\emph{Selecta Math.}
{\bf 6} (2000), 399--442.





\bibitem{HachThesis}
M. Hachtroudi,
\emph{Les espaces d'\'el\'ements \`a connexion projective normale},
Th\`ese de Doctorat \`es Sciences Math\'ematiques,
Facult\'e des Sciences de l'Universit\'e de Paris,
soutenue en juillet 1937,
92 pp.

\bibitem{Hach37}
M. Hachtroudi,
\emph{Les espaces d'\'el\'ements \`a connexion projective normale},
Actualit\'es Scientifiques et Industrielles,
No.~565,
Hermann,
Paris, 1937.

\bibitem{Hach45}
M. Hachtroudi,
\emph{Les espaces normaux},
Publications de l'Universit\'e de T\'eh\'eran,
Tehran, 1945.

\bibitem{Hach48}
M. Hachtroudi,
\emph{Les connexions normales, affines et weyliennes},
Publications de l'Universit\'e de T\'eh\'eran,
Tehran, 1948.

\bibitem{Hach56}
M. Hachtroudi,
\emph{Sur les espaces de Riemann, de Weyl et de Schouten},
Publications de l'Universit\'e de T\'eh\'eran,
No.~297,
Imprimerie de l'Universit\'e de T\'eh\'eran,
Tehran, 1956.

\bibitem{Hach61}
M. Hachtroudi,
\emph{D\=ane\v{s} o Honar}
[Science and Art],
Tehran, 1961.

\bibitem{Hach66}
M. Hachtroudi,
\emph{Tamrinh\=a-ye Riy\=a\.ziy\=at-e Moqaddam\=ati; 
Hendese-ye Dav\=ayer dar Safhe} [Exercises in Elementary Mathematics; 
The Geometry of Circles in the Plane],
Tehran, 1966.

\bibitem{Hach70}
M. Hachtroudi,
\emph{Nazariye-ye A'd\=ad}
[Theory of Numbers],
Tehran, 1970.

\bibitem{Hach83}
M. Hachtroudi,
\emph{Seyr-e Andi\v{s}e-ye Ba\v{s}ar} [The Course of Human Thought], Third Edition, 
Tehran, 1983.

\bibitem{HachAuto}
M. Hachtroudi,
Autobiographical account,
\emph{Majalla-ye Yek\=an},
No.~7,
Tehran,
August 1964.




\bibitem{IAS}
Institute for Advanced Study, ``Mohsen Hachtroudi,''
\textit{School of Mathematics Scholars Archive}.
Available at
\url{https://www.ias.edu/scholars/mohsen-hachtroudi}.


\bibitem{IveyLandsberg}
T. A. Ivey and J. M. Landsberg,
\emph{Cartan for Beginners:
Differential Geometry via Moving Frames and Exterior Differential Systems},
Graduate Studies in Mathematics,
vol.~61,
American Mathematical Society,
Providence, RI, 2003.


\bibitem{LieEngel}
S. Lie,
\emph{Theorie der Transformationsgruppen},
unter Mitwirkung von F. Engel,
3 vols.,
B. G. Teubner,
Leipzig, 1888--1893.

\bibitem{MerkerLie}
J. Merker,
Lie symmetries of partial differential equations and CR geometry,
\emph{J. Math. Sci.}
{\bf 154} (2008), 817--922.

\bibitem{MerkerSabzevari12}
J. Merker and M. Sabzevari,
Explicit expression of Cartan's connection for
Levi-nondegenerate $3$-manifolds in complex surfaces,
and identification of the Heisenberg sphere,
\emph{Cent. Eur. J. Math.}
{\bf 10} (2012), 1801--1835.

\bibitem{MerkerSabzevari14}
J. Merker and M. Sabzevari,
The Cartan equivalence problem for
Levi-nondegenerate real hypersurfaces,
\emph{Izv. Math.}
{\bf 78} (2014), 1158--1194.

\bibitem{MerkerHach}
J. Merker,
Vanishing Hachtroudi curvature and local equivalence to the
Heisenberg pseudosphere,
\emph{Bull. Iranian Math. Soc.}
{\bf 47} (2021), 1775--1792.


\bibitem{MilaniHashtrodi}
A. Milani,
\emph{Eminent Persians: The Men and Women Who Made Modern Iran,
1941--1979},
2 vols.,
Syracuse University Press,
Syracuse, NY, 2008;
see ``Mohsen Hashtrudi,'' vol.~II, pp.~925--929.


\bibitem{Olver}
P. J. Olver,
\emph{Equivalence, Invariants, and Symmetry},
Cambridge University Press,
Cambridge, 1995.

\bibitem{Segre31}
B. Segre,
Intorno al problema di Poincar\'e della rappresentazione
pseudoconforme,
\emph{Rend. Accad. Naz. Lincei}
{\bf 13} (1931), 676--683.

\bibitem{Sharpe}
R. W. Sharpe,
\emph{Differential Geometry:
Cartan's Generalization of Klein's Erlangen Program},
Graduate Texts in Mathematics,
vol.~166,
Springer,
New York, 1997.


\bibitem{Yad1}
H. Soudbakhsh (ed.),
\emph{Y\=adn\=ameh-ye Mohsen Hachtroudi}
(\emph{Memorial Volume for Mohsen Hachtroudi}),
3rd ed.,
Mo'assese-ye Farhang\={\i}-ye Doktor Hachtroudi,
Tehran, 2007. 




\bibitem{Sukhov01}
A. Sukhov,
Segre varieties and Lie symmetries,
\emph{Math. Z.}
{\bf 238} (2001), 483--492.




\bibitem{IranicaHach}
A. S. Tahvildar-Zadeh and F. Majidi,
Ha\v{s}trudi, Mohsen,
\emph{Encyclopaedia Iranica},
Vol.~XII, Fasc.~1 (2003), 52--54.





\bibitem{Tanaka76}
N. Tanaka,
On non-degenerate real hypersurfaces,
graded Lie algebras and Cartan connections,
\emph{Japan. J. Math.}
{\bf 2} (1976), 131--190.


\bibitem{Tresse}
A. Tresse,
Sur les invariants diff\'erentiels des groupes continus
de transformations,
\emph{Acta Math.}
{\bf 18} (1894), 1--88.

\bibitem{VeblenThomas23}
O. Veblen and T. Y. Thomas,
The geometry of paths,
\emph{Trans. Amer. Math. Soc.}
{\bf 25} (1923), 551--608.





\bibitem{Yekan76}
\emph{Majalla-ye Yek\=an},
Special issue devoted to Mohsen Hachtroudi,
vol.~12,
no.~6,
Tehran,
December 1976.

\end{thebibliography}
\end{document}